\documentclass[american,parskip=half,12pt]{scrartcl}

\usepackage[T1]{fontenc}
\usepackage[utf8]{inputenc}
\usepackage{babel}
\usepackage{csquotes}
\usepackage[shortlabels]{enumitem}

\usepackage{amssymb,amsmath,bm,mathtools}
\usepackage{mathabx} 
\usepackage[hyperref,amsmath,thmmarks]{ntheorem}

\usepackage{graphicx}
\usepackage[justification=raggedright]{subcaption}
\usepackage{xcolor}
\usepackage{pgfplots}
\graphicspath{{../images/},{../matlab/images/}}
\usepgfplotslibrary{external}

\usepackage[backend=biber, style=numeric-comp, maxbibnames=9, giveninits=true, backref=false
  ]{biblatex}

\newrobustcmd{\MakeTitleCase}[1]{%
  \ifthenelse{\ifcurrentfield{title}}%
  {\MakeSentenceCase{#1}}%
  {#1}}
\DeclareFieldFormat[article,inbook,incollection,inproceedings]{titlecase}{\MakeTitleCase{#1}}
\AtEveryBibitem{\clearfield{url}}
\AtEveryBibitem{\clearfield{month}}
\renewbibmacro{in:}{%
  \ifentrytype{article}{}{\printtext{\bibstring{in}\printunit{\intitlepunct}}}}

\usepackage[english,pdfusetitle,colorlinks=true,linkcolor=blue,citecolor=green!50!black,urlcolor=blue,bookmarks=true]{hyperref}

\usepackage[ruled]{algorithm2e}

\usepackage{todonotes}
\usepackage{aliascnt}
\usepackage{siunitx}

\newaliascnt{proposition}{lemma}
\newtheorem{proposition}[proposition]{Proposition}
\aliascntresetthe{proposition}

\newaliascnt{corollary}{lemma}

\aliascntresetthe{corollary}

\newaliascnt{theorem}{lemma}
\newtheorem{theorem}[theorem]{Theorem}
\aliascntresetthe{theorem}

\theorembodyfont{\normalfont}
\newaliascnt{definition}{lemma}

\aliascntresetthe{definition}

\newaliascnt{assumption}{lemma}

\aliascntresetthe{assumption}

\newaliascnt{notation}{lemma}

\aliascntresetthe{notation}

\newaliascnt{example}{lemma}

\aliascntresetthe{example}

\newaliascnt{experiment}{lemma}

\aliascntresetthe{experiment}

\newaliascnt{remark}{lemma}

\aliascntresetthe{remark}

\theoremstyle{nonumberplain}
\theoremseparator{:}
\theoremheaderfont{\normalfont\itshape}

\theoremsymbol{\ensuremath{\square}}
\newtheorem{proof}{Proof}

\newcommand{\N}{\ensuremath{\mathbb{N}}}

\newcommand{\R}{\ensuremath{\mathbb{R}}}
\newcommand{\C}{\ensuremath{\mathbb{C}}}

\newcommand{\abs}[1]{\ensuremath{\left\vert#1\right\vert}}

\newcommand{\e}{\mathrm{e}}

\renewcommand{\i}{\mathrm{i}}
\newcommand{\zb}[1]{\ensuremath{\bm{#1}}}

\DeclareMathOperator*{\Div}{div}
\DeclareMathOperator*{\argmin}{arg\,min}
\DeclareMathOperator*{\minimize}{minimize}

\DeclareMathOperator*{\supp}{supp}
\DeclareMathOperator*{\sgn}{sgn}
\DeclareMathOperator{\proj}{proj}
\DeclareMathOperator{\prox}{prox}
\DeclareMathOperator{\grad}{grad}
\DeclareMathOperator{\re}{Re}
\renewcommand{\d}{\, \mathrm{d}}

\newcommand{\norm}[1]{\left\lVert #1
  \right\rVert}

\newcommand{\dd}{\, \mathrm{d} }
\newcommand{\bn}{{\bm{n}}}
\newcommand{\bc}{{\bm c}}

\newcommand{\bl}{{\bm \ell}}
\newcommand{\bg}{{\bm g}}
\newcommand{\br}{{\bm r}}

\newcommand{\bx}{{\bm x}}
\newcommand{\by}{{\bm y}}
\newcommand{\bk}{{\bm k}}

\newcommand{\bv}{{\bm v}}
\newcommand{\bw}{{\bm w}}

\newcommand{\be}{{\bm e}}
\newcommand{\bz}{{\bm z}}
\newcommand{\bo}{{\bm 0}}
\newcommand{\bh}{{\bm h}}

\newcommand{\ui}{u^{\mathrm{inc}}}

\newcommand{\utot}{u^{\mathrm{tot}}}

\newcommand{\us}{u^{\mathrm{sca}}}
\newcommand{\vs}{\bm{v}^{\mathrm{sca}}}

\newcommand{\rs}{r_{\mathrm{s}}}

\newcommand{\rM}{r_{\mathrm{M}}}
\newcommand{\lM}{L_{\mathrm{M}}}
\newcommand{\ls}{L_{\mathrm{s}}}
\newcommand{\dft}{\bm{F}_{\mathrm{DFT}}}
\newcommand{\ndft}{\bm{F}_{\mathrm{NDFT}}}

\title{Total variation-based reconstruction and phase retrieval for diffraction tomography}

\date{}
\author{%
  Robert Beinert\thanks{TU Berlin, Institute of
    Mathematics, MA 4-3, Straße des 17. Juni 136, D-10623 Berlin, Germany.
  }\\%
  {\footnotesize\href{mailto:beinert@math.tu-berlin.de}{beinert@math.tu-berlin.de}}
  \and 
  Michael Quellmalz\footnotemark[1]\\%
  {\footnotesize\href{mailto:quellmalz@math.tu-berlin.de}{quellmalz@math.tu-berlin.de}}
}

\begin{document}

\def\sectionautorefname{Section}
\def\subsectionautorefname{Section} 

\maketitle

\begin{abstract}
  In optical diffraction tomography (ODT), the three-dimensional
  scattering potential of a microscopic object rotating around its
  center is recovered by a series of illuminations with coherent
  light.  Reconstruction algorithms such as the filtered
  backpropagation require knowledge of the complex-valued wave at the
  measurement plane, whereas often only intensities, i.e., phaseless
  measurements, are available in practice.

  We propose a new reconstruction approach for ODT with unknown phase
  information based on three key ingredients.  First, the light
  propagation is modeled using Born's approximation enabling us to use
  the Fourier diffraction theorem.  Second, we stabilize the inversion
  of the non-uniform discrete Fourier transform via total variation
  regularization utilizing a primal-dual iteration, which also yields
  a novel numerical inversion formula for ODT with known phase. The
  third ingredient is a hybrid input-output scheme.  We achieved
  convincing numerical results, which indicate that ODT with phaseless
  data is possible.
  The so-obtained 2D and 3D reconstructions are even comparable to the ones
  with known phase.
  
  \medskip
  \noindent
  \textit{Math Subject Classifications.}
  42B05,  
  47J06, 
  65T50,  
  92C55 

\end{abstract}

\section{Introduction}

Illuminating three-dimensional (3D) objects from different angles, we are
interested in tomographic reconstructions.  In classical X-ray
computerized tomography, the light is assumed to propagate along
straight lines enabling reconstruction via the inverse Radon
transform.  This assumption relies on the fact that the employed
X-rays have a very short wavelength, which is much smaller than any
object features of interest.  In \emph{optical diffraction tomography}
(ODT), the image is obtained via visible light of wavelengths of
hundreds of nanometers being the size of studied features such as
subcellular structures.  Therefore, we rely on a more sophisticated
light propagation model taking the occurring diffraction into account.

Mathematically speaking, we aim to recover the three-dimensional
scattering potential of the unknown object given two-dimensional
measurements of the complex wave function at a plane behind the object
for different illumination directions. 
We assume that the illumination directions are known beforehand,
motion detection methods were discussed in \cite{ElbRitSchSchm20}. To model the light
propagation, we take advantage of the Born approximation, which is valid
for objects of small size and variation of the scattering potential.
A well-established inversion method for the Born-based ODT model is
the backpropagation \cite{Dev82}, which is still widely used in
practice \cite{ODTbrain,FaSmLiSa17}.  The backpropagation, however,
relies on complex-valued measurements whereas only phaseless
intensities are available in most practical setups.  
This means that the object reconstruction becomes more ill posed and
is coupled with a phase retrieval problem.

Phase retrievals of different kinds have been intensely studied in the
literature.  In general, solutions are highly ambiguous.
Characterizations of the whole solution set and uniqueness guarantees
under additional constraints like known support or non-negativity have
been studied for functions with continuous domains
\cite{Aku56,KliSacTik95, KliKam14,Wal63,GroKopRat20}, where the
measurements correspond to the intensity of the Fourier transform.
Analogous results exist for fractional Fourier intensities
\cite{Jam14} for instance.  Even if the Fourier phase retrieval
problem is discretized, the ambiguousness remains an issue
independently of how many measurements are available
\cite{BenBeiEld17, BeiPlo15, GroKopRat20, BarEpsGreMag20,
  Edi19,BeiPlo20}.  To overcome these problems, more general phase
retrieval problems, where the intensity of arbitrary linear
measurements are given, have been studied.  Under certain assumptions
on the linear measurements, the phase retrieval problem becomes unique
\cite{CanStrVor13, CanEldStrVor13, FanLia13, FanStr20}.  Using a
tensorial lifting, the problem may be exemplarily solved using
semi-definite programming \cite{CanStrVor13, CanEldStrVor13},
tensor-free primal-dual methods \cite{BeiBre21}, or polarization
techniques \cite{AleBanFicMix14, PfaSal19}.  Moreover, the phase
retrieval problem is well understood in the context of frames
\cite{GroKopRat20, PfaSal19, AleBanFicMix14}.

Unfortunately, phase retrieval in ODT does not fit in the above
setups since the given measurements with respect to the Born
approximation are non-linear.  Therefore, the well-established
theoretical results regarding ambiguousness, uniqueness, and stability
cannot be transferred to our problem.  In X-ray tomography,
where the intensity measurements are modeled using the Radon transform
of the unknown object followed by diffraction, uniqueness of the
reconstruction is guaranteed \cite{Mar18}.  The difference in the
employed mathematical model, however, prevents a direct application of
the uniqueness results in ODT.

For retrieving an unknown phase in ODT, different approaches have been
studied.  An iterative method of propagation and backpropagation in
free space between the measurements plane and a parallel plane through
the object was proposed by Maleki and Devaney \cite{MalDev93}.  The
Gabor holography compared unfavorably in numerical tests
\cite{WedSta95}.  A phase retrieval method utilizing ideas of Fourier
ptychography was developed in \cite{HorChuOuZheYan16}.  Other
approaches rely on the acquisition of more data using two parallel
measurement planes \cite{GbuWol02} or the phase shifting
interferometry (PSI) requiring 4 intensity images per illumination
direction.  Furthermore, phase retrieval based on a far zone
approximation of Born \cite{CheHan01} or using the X-ray transform
\cite{GurDavPogMayWil04} has been studied.

In this paper, we consider the phase retrieval problem in ODT, where
the absolute value of the wave is measured at one plane per
illumination direction, and we assume that the Born approximation is
valid.  Our proposed method is based on the hybrid input-output
algorithm from Fourier phase retrieval combined with a total variation
(TV) regularization utilizing a primal-dual method.  TV
regularization has already been applied to limited angle
reconstruction in ODT utilizing an approximate Radon transform
\cite{Kuj14}.

The main contributions of this paper are as follows:
\begin{itemize}
\item We consider the inverse nonuniform discrete Fourier transform
  (NDFT) as an inverse problem and apply TV-based regularization
  methods to improve the inversion especially in the presence of
  noise.
  We apply this regularized NDFT inversion to obtain a novel
  reconstruction method for ODT with known phase.
\item We propose a numerical algorithm for solving the phase retrieval
  problem in ODT.  This algorithm makes use of an adaption of the
  hybrid input-output scheme and the TV-regularized inversion of the
  NDFT.
\item Our numerical simulations
  indicate that phase retrieval is possible and greatly benefits from
  the regularized ODT inversion.  
  In total, the reconstruction quality
  is comparable to the known-phase inversion.
\end{itemize}

\textbf{Outline of the paper.}  The employed notation is introduced in
\autoref{sec:preliminary-remarks}.  The basics on diffraction
tomography and the theoretical foundation of the inversion methods are
explained in \autoref{sec:odt}.  We introduce the discretized setting
used for implementations in \autoref{sec:discretization}.  Since the
inversion of the NDFT is the key element for reconstruction in
diffraction tomography, we discuss different computation methods for
the inverse NDFT in \autoref{sec:indft}.  On the basis of these
methods, different numerical ODT inversions for complex-valued data
are studied and compared.  The proposed phase retrieval method for ODT
is derived in \autoref{sec:phase-retrieval-odt}.  Finally, numerical
simulations in \autoref{sec:numer-pr} show that phase retrieval in ODT
is possible and we achieve favorable reconstruction results.
Conclusions are given in
\autoref{sec:conclusion}.

\section{Notation and preliminary remarks}
\label{sec:preliminary-remarks}

In the following, we use the (symmetric) index set
\begin{equation*}
  \mathcal I_K \coloneqq \{-K/2, -K/2+1,\dots,K/2-1\},
\end{equation*}
where $K\in2\N$.  Although $\mathcal I_K$ may be adapted for odd $K$,
we only consider the even case for convenience.  The set of the first
$M\in\N$ positive integers is further denoted by 
\begin{equation*}
  \mathcal J_M \coloneqq \{1, 2,\dots,M\}.
\end{equation*}
Numerically, we represent the unknown object as \emph{$d$-dimensional
  signal}
\begin{equation*}
  \bm x \coloneqq [x_{\bm k}]_{\bm k \in \mathcal I_K^d} \in
  \mathbb R^{K^d},
\end{equation*}
which may be interpreted as samples of a multivariate function.  The
\emph{(weighted) 2-norm} of the $d$-dimensional signal $\bx$ is
defined by
\begin{equation*}
  \| \bm x \|_2^2
  \coloneqq
  \sum_{\bm k \in \mathcal I_K^d} | x_{\bm k}|^2
  \qquad\text{and}\qquad
  \| \bm x \|_{2,\bm w}^2
  \coloneqq
  \sum_{\bm k \in \mathcal I_K^d} w_{\bm k} | x_{\bm k}|^2
\end{equation*}
for some non-negative weights $\bm w \in \mathbb R_{\ge 0}^{K^d}$.

Besides the signals, we also require the notion of a $d$-dimensional
\emph{discrete vector field} $\bm y \in \mathbb R^{K^d \times d}$, where the
elements of $\bm y$ are itself $d$-dimensional vectors, i.e.\
$y_{\bm k} \in \mathbb R^d$ for $\bm k \in \mathcal I_K^d$.  The
signal $\bm y_{\cdot,\ell}$ consisting of the $\ell$th
components is defined by
$[\bm y_{\cdot,\ell}]_{\bm k} \coloneqq [y_{\bm k}]_\ell$ for $\bm k \in
\mathcal I_K^d$ and $\ell \in \mathcal J_d$.  The \emph{(1,2)- and (2,2)-norm} of
a discrete vector field is given by
\begin{equation*}
  \| \bm y \|_{1,2} \coloneqq \sum_{\bm k \in \mathcal I_K^d}
  \|y_{\bm k} \|_2
  \qquad\text{and}\qquad
  \| \bm y \|_{2,2}^2 \coloneqq \sum_{\bm k \in \mathcal I_K^d}
  \|y_{\bm k} \|_2^2.
\end{equation*}
We use analogous definitions for signals whose domains are not
$\mathcal I_K^d$.

Finally, we interpret all occurring complex vector spaces
$\mathbb C^N$ as $2N$-dimensional vector spaces over $\R$ with the
inner product
\begin{equation*}
   \langle \bx,
   \by \rangle \coloneqq \sum_{n=1}^N \re [x_n \bar y_n ]
\end{equation*}
for $\bm x, \bm y \in \mathbb C^N$, which is necessary to apply
certain techniques and numerical methods from convex analysis.  For
$\bm A \in \mathbb C^{N \times K}$ declaring the mapping
$\bm x \mapsto \bm A \bm x$ from $\mathbb R^K$ to $\mathbb C^N$, the
adjoint mapping becomes $\by \mapsto \re[\bm A^* \by]$ from $\mathbb
C^N$ to $\mathbb R^K$.

\section{Diffraction tomography}
\label{sec:odt}

In this section, we describe the setup in \emph{optical diffraction
  tomography} (ODT), cf.\ \cite{KirQueRitSchSet21}.  We want to
reconstruct an object's \emph{scattering potential}
$f\colon\R^d \to \R$, where we restrict our considerations to the
practically meaningful dimensions $d\in\{2,3\}$.  The object is
embedded in a homogeneous background and is illuminated by a
monochromatic plane wave.  The
\emph{incident field}
\begin{equation} \label{eq:uinc}
  \ui(\bx)=\e^{\i k_0 x_d}, \qquad\bx\in\R^d,
\end{equation}
with \emph{wave number} $k_0>0$ propagates in the direction of $x_d$.
The scattering potential admits the relation
$
f(\bx)
=
k_0^2\, ({n(\bx)}/{n_0} -1 ),
$
where $n(\bx)$ is the space-dependent refractive index and $n_0$ is the constant refractive index of the surrounding medium.
In most practical cases, such as the imaging of biological cells in water,
the refractive index in the object is greater than in its surrounding,
cf.\ \cite{Liu16},
so we will assume that $f\ge0$.

For the tomographic reconstruction, we rotate the unknown object
during the measurement process.  The rotation at time $t\in [0,T]$ for
$T>0$ is described by the known rotation matrix
$ R_t\in \mathrm{SO}(d) $.  The rotated object has the scattering
potential~$f(R_t \cdot)$. The resulting total field at time $t$ is
denoted by $\utot_t\colon\R^d\to\C$.  We decompose the total field
into
$$
\utot_t
=
\ui + \us_t
$$
with the known incident field $\ui$ given in \eqref{eq:uinc} and the
scattered field $\us_t$.

\begin{figure}[t]
  \begin{center}
    \begin{tikzpicture}[scale=0.5]
      
      \draw[->] (-5,0)--(5,0) node[right] {$x_1$};
      \draw[->] (0,-4.5)--(0,4.5) node[right] {$x_3$};
      \draw[->,opacity=.7] (-2,-2)--(2,2) node[right] {$x_2$};
      
      \fill[red,opacity=0.3] (0,0) circle (2);
      \node at (1,-0.8) {$f$};
      \draw[->]  (-2.46, 0.43) arc (170:100:2.5);
      
      \draw[dashed] (-5,-3.5) -- (5,-3.5);
      \draw[dashed] (-5,-4) -- (5,-4);
      \draw[dashed] (-5,-4.5) -- (5,-4.5);
      \node at (-8,-4) {\footnotesize{incident field}};
      \draw[->] (5.5,-4.5) -- (5.5,-3.5) node[right] {$\ui$};
      
      \draw[line width = 3pt] (5,3.5) -- (-5,3.5) node [left] {\footnotesize{measurement plane}};
      \node at (7,3.5) {$x_3 = \rM$};
    \end{tikzpicture}
  \end{center}
  \caption{Experimental setup for $d=3$ with measurement plane located at $x_3=\rM$.
  }	\label{fig:trans}
\end{figure}
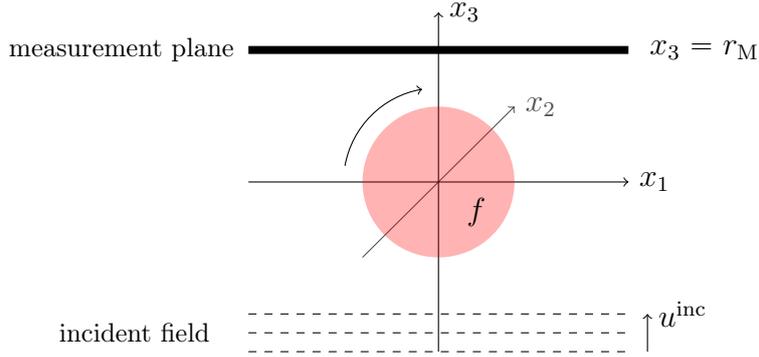

The total field $\utot_t$ at time $t$ is measured at a fixed plane
behind the object.  More precisely, we detect either
$\utot_t(\bx',\rM)$ or $|\utot_t(\bx',\rM)|$ for all $\bx'\in\R^{d-1}$
and a fixed $\rM>0$; the measurement plane is thus orthogonal to the $x_d$ axis.
Since the plane lies outside
the object, $f$ is assumed to be bounded and compactly supported in
the ball $\mathcal B_{\rM}^d = \{\bx\in\R^d: \norm{\bx}_2<\rM\}$.  The
experimental setup is illustrated in \autoref{fig:trans}.

\subsection{Born approximation}

To model the scattered wave $\us_t$, we employ the first order
\emph{Born approximation} \cite{KakSla}, where $\us_t$ solves the
\emph{Helmholtz equation}
\begin{equation} \label{eq:Helmholtz}
  (\Delta + k_0^2)\, \us_t
  = - f \ui.
\end{equation}
Under the Sommerfeld radiation condition, which basically says that
the scattered field is an outgoing wave, the Helmholtz equation
\eqref{eq:Helmholtz} has a unique solution, see \cite[Chap. 2]{ColKre13}.  This solution may be
computed by using the outgoing \emph{Green function} $G_d$ of the
Helmholtz equation, see \cite[Sect.\
3.3]{NaWue00}, which is given by
\begin{equation*}
  G_3(\bx)
  \coloneqq
  \frac{\e^{\i k_0 \norm{\bx}_2}}{4\pi \norm{\bx}_2}
  ,\qquad \bx\in\R^3\setminus\{\bo\},
\end{equation*}
and
\begin{equation*}
  G_2(\bx)
  \coloneqq \frac{\i}{4} H^{(1)}_0(k_0 \norm{\bx}_2)
  ,\qquad \bx\in\R^2\setminus\{\bo\},
\end{equation*}
where 
$H^{(1)}_0(x)
\coloneqq
(\pi\i)^{-1} 
\int_{-\infty}^{\infty}
\e^{\i x \cosh t} \dd t
$, 
$x>0$, 
denotes the zeroth order Hankel function of the
first kind \cite[\nopp(8.421.1)]{GrRy07}.
The scattered field $\us_t$ is then given by the convolution
\begin{equation} \label{eq:conv}
  \us_t(\bx)
  = \left[(f(R_t\cdot)\, \ui) * G_d\right] (\bx)
  \coloneqq \int_{\R^d} f(R_t \by)\, \ui(\by)\, G_d(\bx-\by)\d\by
  ,\qquad \bx\in\R^d.
\end{equation}
Restricting the convolution to the measurement plane, we model the
given ODT data using the (non-linear) forward map
$\mathcal D^{\mathrm{tot}}$ defined by
\begin{equation} \label{eq:Dtot-op}
  \mathcal D^{\mathrm{tot}} f(\bx',t) 
  \coloneqq [(f(R_t \cdot)\, \ui) * G_d](\bx', \rM)
  + \ui(\bx',\rM),
  \qquad \bx'\in\R^{d-1},\ t\in[0,T].
\end{equation}

\subsection{Fourier diffraction theorem}

Besides the convolutional representation \eqref{eq:Dtot-op}, the
forward map $\mathcal D^{\mathrm{tot}}$ may be represented using the
Fourier transform, which significantly simplifies the inversion.
The $d$-di\-men\-sion\-al \emph{Fourier transform} is defined
by
\begin{equation} \label{eq:F}
  \mathcal F [f](\by)
  \coloneqq (2\pi)^{-d/2} \int_{\R^d} f(\bx)\,\e^{-\i\, \bx\cdot\by} \dd \by
  ,\qquad \by\in\R^d,
\end{equation} 
and the \emph{partial Fourier transform} with respect to the first $d-1$
coordinates by
\begin{equation} \label{eq:F1}
\mathcal F_{1,\dots,d-1} [f](\by',x_d)
\coloneqq (2\pi)^{-(d-1)/2} \int_{\R^{d-1}} f(\bx',x_d)\, \e^{-\i\, \bx'\cdot\by'} \dd \bx'
\end{equation}
for $(\by',x_d)\in\R^{d-1}\times\R$.  Recalling that $k_0$ is the
wave number of the incident field, we further define
\begin{equation*}
\kappa(\by') \coloneqq \sqrt{k_0^2-\norm{\by'}_2^2}
\qquad\text{and}\qquad
\bh(\by')
\coloneqq
\begin{pmatrix}
  \by'\\\kappa(\by')-k_0
\end{pmatrix}
\in\R^d
.
\end{equation*}
where $\by'\in \mathcal B_{k_0}^{d-1}$.  The \emph{Fourier diffraction
  theorem} \cite{KakSla,NaWue00} now states
\begin{equation} \label{eq:recon}
  \mathcal F_{1,\dots,d-1}[ \us_t](\by', \rM)
  = 
  \frac{\i\,\sqrt{\pi}}{\sqrt2\, \kappa(\by')}\,
    \e^{\i\, \kappa(\by')\, \rM}\,
  \mathcal F [f](R_t\bh(\by'))
  ,\qquad \by' \in \mathcal B^{d-1}_{k_0},\ t\in [0,T].
\end{equation}
In general, the identity \eqref{eq:recon} only holds in a
distributional sense \cite{KirQueRitSchSet21}.  Geometrically, the set
$\{\bh(\by') : \norm{\by'}_2<k_0\}$ is a semisphere or semicircle,
whose center is on the negative $x_d$ axis, and which passes through
the origin. On the right-hand side of \eqref{eq:recon}, we thus
evaluate the Fourier transform $\mathcal F[f]$ on a rotated semisphere
or semicircle; on the left-hand side, partial Fourier transforms of
the measurements $\us_t(\cdot,\rM)$ are computed.

\subsection{Uniqueness of the reconstruction}

Under certain conditions, the scattering potential $f$ is uniquely
determined by the given measurements of $\utot_t$.  The following
uniqueness result is based on the Paley--Wiener theorem \cite[Thm.\
7.3.1]{Ho90}.

\begin{proposition}[Paley--Wiener] \label{thm:pw} Let $f\in L^1(\R^d)$
  have compact support.  Then the Fourier transform
  $\mathcal F[f]\in C(\R^d)$ is uniquely extendable to an analytic
  function on $\C^d$.
\end{proposition}

\begin{theorem}[ODT inversion] \label{thm:unique}
Let $f\in L^1(\R^d)$ have a compact support in $\mathcal B_{\rM}^d$,
and let the  Lebesgue measure of
\begin{equation} \label{eq:Y}
  \mathcal Y 
  \coloneqq 
  \{ R_t \bh(\by') : \by'\in\mathcal B_{k_0}^{d-1},\, t\in [0,T]\}
  \subset \R^d
\end{equation}
be positive.  Then $f$ is uniquely determined by
$\mathcal D^{\mathrm{tot}} [f]$ in \eqref{eq:Dtot-op}.
\end{theorem}

\begin{proof}
  Given $\utot_t$, we have $\us_t = \utot_t - \e^{\i k_0 \rM}$.  By
  the Fourier diffraction theorem \eqref{eq:recon}, we obtain
  $\mathcal F[f]$ on $\mathcal Y$, which has positive Lebesgue
  measure.  Any $d$-variate analytic function is uniquely determined
  by its values on a set of positive Lebesgue measure \cite[Sect.\
  4.1]{KraPar02}.  Since $\mathcal F[f]$ is an analytic function by
  \autoref{thm:pw}, $\mathcal F[f]$ is uniquely determined on $\R^d$
  and thus $f$ itself due to the uniqueness of the Fourier transform
  \cite[Thm.\ 4.47]{PlPoStTa18}.
\end{proof}

The requirements of \autoref{thm:unique} are, for instance, fulfilled
if $R_t$ describes a full rotation around an axis different from $x_3$
in 3D, or if the rotation angle covers an interval in 2D.

\section{Discretization}
\label{sec:discretization}

The mapping $\mathcal D^{\mathrm{tot}}$ and its inversion can be
numerically realized utilizing the Fourier diffraction theorem
\eqref{eq:recon}.  For this purpose, we first restrict the total field
$u^{\mathrm{tot}}_t$ or, more precisely, the scattered field
$u^{\mathrm{sca}}_t$ to the cube
$[-\lM,\lM]^{d-1}$
and then discretize the
fields using $N\in 2\N$ samples for each spatial direction and
$M \in \N$ time steps, cf.\ \cite{KirQueRitSchSet21}.  On the uniforms grids
\begin{equation} \label{eq:z'}
\bz'_\bn 
\coloneqq \tfrac{2\lM}{N} \, \bn 
,\quad \bn\in\mathcal I_N^{d-1}
\qquad\text{and}\qquad
t_m \coloneqq \tfrac{T}{M} \, m ,
\quad m\in\mathcal J_M,
\end{equation}
the scattered fields $u_t^{\mathrm{sca}}(\bz',\rM)$ are represented by the vectors
\begin{equation*}
  \vs_m
  \coloneqq \left[v_{m,\zb n}^{\mathrm{sca}}\right]_{\zb n\in\mathcal I_N^{d-1}}
  \quad\text{ with }\quad
  \vs_{m,\zb n}
  \coloneqq
  \us_{t_m} (\bz_{\bn}',\rM ).
\end{equation*}
The $(d-1)$-dimensional \emph{discrete Fourier transform} (DFT) of
$\bm v_m^{\mathrm{sca}}$ is given by
\begin{equation} \label{eq:dft}
  \left[\dft \, \vs_m\right]_\bl
  \coloneqq
  \sum_{\zb n \in\mathcal I_N^{d-1}}
  v_{m,\zb n}^{\mathrm{sca}} \, \e^{-2\pi\i \zb n\cdot\bl/N}
  ,\qquad \bl\in\mathcal I_N^{d-1},
\end{equation}
and may be inverted using the inverse DFT
\begin{equation} \label{eq:idft}
  v_{m,\zb n}^{\mathrm{sca}}
  =
  \frac{1}{N^{d-1}}
  \sum_{\bl\in\mathcal I_N^{d-1}}
  \left[\dft \, \bm v_m^{\mathrm{sca}} \right]_\bl
  \e^{2\pi\i \zb n\cdot\bl/N}
  ,\qquad \zb n \in\mathcal I_N^{d-1}.
\end{equation}
The left-hand side of the Fourier diffraction theorem \eqref{eq:recon}
may thus be discretized as
\begin{equation} \label{eq:F1-dft}
\mathcal F_{1,\dots,d-1} \us_{t_m}\left(\by'_\bl,\rM\right)
\approx
(2\pi)^{-\frac{d-1}{2}} \left(\tfrac{2\lM}{N}\right)^{d-1}
[\dft \, \zb v^{\mathrm{sca}}_{m}]_{\bl}
,\qquad \bl\in\mathcal I_N^{d-1},
\end{equation}
where $\by'_\bl \coloneqq \frac{\pi}{\lM} \, \bl$ with
$\bl\in\mathcal I_N^{d-1}$, which provides an approximation of
$\mathcal F[f]$ on the nonuniform grid
\begin{equation} \label{eq:Y_MN}
  \mathcal Y_{M,N}
  \coloneqq
  \left\{
  R_{t_m} \bh(\by'_\bl)
  :
  m\in \mathcal J_M, \bl\in\mathcal I_N^{d-1},
  \norm{\by_{\bl}'} < {k_0} \right\}.
\end{equation}
For simplicity, we assume $\norm{\by'_\bl}_2 < k_0$ for all
$\bl\in\mathcal I_N^{d-1}$ and truncate the grid to the
points satisfying this assumption.

For the right-hand side of the Fourier diffraction theorem, we assume
$\supp f \subset [-\ls,\ls]^d$ and use $K \in 2\N$ sampling points per
direction.  Using the uniform grid
$$
\bx_\bk \coloneqq \tfrac{2\ls}{K} \, \bk
,\qquad \bk\in\mathcal I_K^d,
$$ 
we approximate the scattering potential $f$ by
$$
\bm f 
= \left[f_{\bk}\right]_{\bk\in\mathcal I_K^d}
\coloneqq \left[f(\bx_\bk)\right]_{\bk\in\mathcal I_K^d}\in\R^{K^d}.
$$
The $d$-dimensional \emph{nonuniform discrete Fourier transform}
(NDFT) of $\bm f$ onto the set $\mathcal Y_{M,N}$ is defined by
\begin{equation} \label{eq:ndft}
  [\ndft \bm f]_{m,\zb \ell}
  \coloneqq
  \sum_{\bk\in\mathcal I_K^{d}}
  f_{\bk}\, \e^{-\i \bx_\bk\cdot \left(R_{t_m}\bh(\by'_\bl)\right)}
  ,\qquad m\in\mathcal J_M,\ \bl\in \mathcal I_{N}^{d-1},
\end{equation}
which approximates the right-hand Fourier transform in
\eqref{eq:recon} via
\begin{equation} \label{eq:F-ndft}
\mathcal F[f]\left(R_{t_m} \bh\left(\by'_\bl\right)\right)
\approx
(2\pi)^{-\frac{d}{2}} \left(\frac{2\ls}{K}\right)^{d}
\left[\ndft \bm f\right]_{m,\zb\ell}
,\qquad m\in\mathcal J_M, \zb\ell\in \mathcal I_{N}^{d-1}.
\end{equation}

Inserting the approximations \eqref{eq:F1-dft} and \eqref{eq:F-ndft}
into the Fourier diffraction theorem \eqref{eq:recon}, we obtain
\begin{equation*}
  (2\pi)^{-\frac{d-1}{2}} \left(\frac{2\lM}{N}\right)^{d-1}
  [\dft \vs_{m}]_{\bl}
  \approx
  \frac{\i\,\sqrt{\pi}}{\sqrt2\, \kappa(\by'_\bl)}\,
  \e^{\i\, \kappa(\by'_\bl)\, \rM}\,
  (2\pi)^{-\frac{d}{2}} \left(\frac{2\ls}{K}\right)^{d}
  \left[\ndft \bm f\right]_{m,\zb\ell}.
\end{equation*}
With the definition
$$
\bc = [c_{\bl}]_{\bl\in\mathcal I_N^{d-1}} \in \R^{N^{d-1}}
\quad\text{ with }\quad
c_{\bl}
\coloneqq
\frac{\i}{\kappa(\by'_\bl)}\,
\e^{\i\, \kappa(\by'_\bl)\, \rM}\,
\left(\frac{N}{\lM}\right)^{d-1}
\left(\frac{\ls}{K}\right)^{d}
,
$$
the discrete Fourier diffraction theorem becomes
\begin{equation*}
[\dft \vs_{m}]_{\bl}
\approx
c_{\bl}\,
\left(\ndft \bm f\right)_{m,\zb\ell}.
\end{equation*}
Incorporating the incident field, we approximate $\mathcal D^{\mathrm{tot}}$ in \eqref{eq:Dtot-op} by the discrete diffraction tomography map
\begin{equation} \label{eq:Dtot}
\bm D^{\mathrm{tot}} \bm f
\coloneqq \dft^{-1} (\bc\odot\ndft \bm f)
+ \e^{\i k_0 \rM},
\qquad\bm f\in \R^{K^d},
\end{equation}
where $\odot$ denotes the Hadamard (entrywise) product of two matrices
or vectors.  The numerical computation of the discrete forward
transform $\bm D^{\mathrm{tot}}$ basically consists of an NDFT, an
entrywise multiplication with $\bc$, and a DFT, see \autoref{alg:fw}.

\begin{algorithm}
  \KwIn{Scattering potential $\bm f = \left[f(\bx_\bk )\right]_{\bk\in\mathcal I_K^d}$.}
  
  Compute the NDFT $\bm g \coloneqq \ndft \bm f$  by \eqref{eq:ndft}\;
  
  \For{$m\in\mathcal J_M$}{
    $\bm{\tilde g}_m \coloneqq \bc \odot \bm g_m$\;
    
    $\bm v_m^{\mathrm{sca}} \coloneqq \dft^{-1} \bm{\tilde g}_m$\;
  }
  
  \KwOut{Scattered wave
    $\utot_{t_m}\!\left(\bn\frac{2\lM}{N},\rM\right) \approx \bv_{m,\bn}^{\mathrm{sca}} + \e^{\i k_0 \rM}$,
    $m\in\mathcal J_M$, $\bn\in\mathcal I_N^{d-1}$.}
  
  \caption{Computation of
    $\bm D^{\mathrm{tot}}$}
  \label{alg:fw}
\end{algorithm}

\begin{figure}[t]\centering
\begin{tikzpicture}[scale=.85]
  \tikzset{font=\small}
  \begin{scope}[shift={(0,0,2)}]
    \draw[->] (0,0,0)--(1,0,0) node[right]{$x_1$};
    \draw[->] (0,0,0)--(0,0,-1) node[right]{$x_2$};
  \end{scope}
  
  \draw[fill=yellow!90!black] (0,1,0) -- (2,1,0)--(2,1,2)--(0,1,2)--cycle;
  \node[below] at (2,.8,0) {$\vs_1$};
  \draw[fill=red!90!black] (0,2.5,0) -- (2,2.5,0)--(2,2.5,2)--(0,2.5,2)--cycle;
  \node[below] at (2,2.3,0) {$\vs_2$};
  
  \draw[->] (2.6,1,0)--(4,1,0) node[above left]{$\dft$} node[below left]{$\oslash\, \bc$};
  
  \begin{scope}[shift={(5,0)}]
  \draw[fill=yellow!60!green,opacity=.9] (0,1,0) -- (2,1,0)--(2,1,2)--(0,1,2)--cycle;
  \node[below] at (2,.8,0) {$\bg_1$};
  \draw[fill=red!70!blue,opacity=.9] (0,2.5,0) -- (2,2.5,0)--(2,2.5,2)--(0,2.5,2)--cycle;
  \node[below] at (2,2.3,0) {$\bg_2$};
  
  \draw[->] (2.6,1)--(3.8,1) node[above left] {$R_t \bh$};
  \end{scope}
  
  \begin{scope}[shift={(9.5,1.5,0)},scale=.8]
    \draw[dashed] (0,-2) arc (-90:0:2cm);	
    \draw (2,0) arc (0:90:2cm);
    \draw (0,2) arc (90:270:.5cm and 2cm);
    \draw[dashed] (0,-2) arc (-90:90:.5cm and 2cm);
    
    \begin{scope}
      \clip(3,3) -- (0,3) -- (0,2) arc (90:270:.5cm and 2cm) -- (0,-3) -- (3,-3) -- (3,3);
      \shade[ball color=red!70!blue, opacity=0.70] (0,0) circle (2cm);
    \end{scope}
  
    \draw (0,-2) arc (180:0:2cm);	
    \draw (0,-2) arc (180:360:2cm and .5cm);
    \draw[dashed] (0,-2) arc (180:0:2cm and .5cm);
    
    \begin{scope}
      \clip(0,0) -- (0,3) -- (0,-2) arc (180:360:2cm and .5cm) -- (4,3) -- (0,3) -- (0,0);
      \shade[ball color=yellow!60!green, opacity=0.70] (2,-2) circle (2cm);
    \end{scope}
  \end{scope}

  \node at (11,-.9,0) {$\mathcal Ff$};
  
  \begin{scope}[shift={(14.5,0)}]
    \draw[->] (-2,1)--(-.3,1) node[above left] {inverse} node[below left]{NDFT};
    
    \shade[yslant=0,xslant=0,left color=blue!30!white,
    right color=blue!80] (0,0) rectangle +(2,2);
    \shade[yslant=0,xslant=1,left color=blue!40!white,
    right color=blue!80] (-2,2) rectangle +(2,1);
    \shade[yslant=1,xslant=0,left color=blue!80,
    right color=blue!70] (2,-2) rectangle +(1,2);
    \draw (0,0) -- (2,0) -- (2,2) -- (0,2)--cycle;
    \draw (0,2) -- (2,2) -- (3,3) -- (1,3)--cycle;
    \draw (2,0) -- (3,1) -- (3,3) ;
    \draw[dashed] (0,0) -- (1,1) -- (3,1) ;
    \draw[dashed] (1,1) -- (1,3) ;
    
    \node at (1,-.9,0) {$\zb f$};
  \end{scope}
  
\end{tikzpicture}

\caption{Schematic overview of the reconstruction process in 3D using \autoref{alg:bw}.
  We exemplarily use two rotations of angle $0$ and $\pi/2$ around the $x_2$ axis.
  The first step consists of 2D DFTs along each layer of the data
  $\vs$ and the division by $\bc$.
   These functions are aligned on semispheres in the Fourier space according to the map $R_t\bh$.
  Finally, $\zb f$ is recovered by a 3D inverse NDFT.
  \label{fig:bw}}
\end{figure}
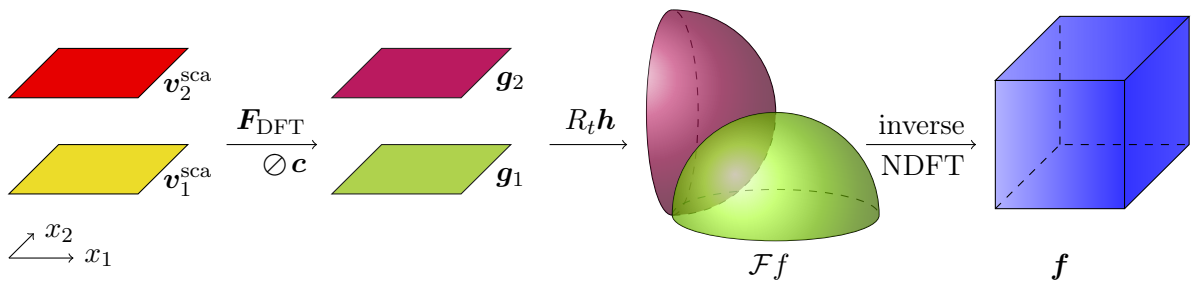

The reconstruction of $\bm f$, see \autoref{alg:bw}, is basically an
inverse of the forward transform.  The first two parts are very simple
to invert utilizing the DFT \eqref{eq:idft} and a pointwise division.
However, the inversion of the NDFT is a more challenging problem that
is covered in the following \autoref{sec:indft}.  A schematic overview
about the reconstruction steps is provided in \autoref{fig:bw}.

\begin{algorithm}
  \KwIn{Scattered wave 
    $\bv_{m,\bn}^{\mathrm{sca}} \coloneqq \utot_{t_m}\!\left(\bn\frac{2\lM}{N},\rM\right) - \e^{\i k_0 \rM}$,
    $m\in\mathcal J_M$, $\bn\in\mathcal I_N^{d-1}$.
  }
  
  \For{$m\in\mathcal J_M$}{
    $\bm{\tilde g}_m \coloneqq \dft \bm v_m^{\mathrm{sca}}$\;
    
    $\bm g_m \coloneqq \bm{\tilde g}_m \oslash \bc $, where $\oslash$ denotes the Hadamard (entrywise) division\; 
  }
  
  Compute the inverse NDFT by solving $\ndft \bm f =\bm g$ for $\bm f$, see \autoref{sec:indft}\;
  
  \KwOut{Scattering potential $\bm f \approx \left[f(\bx_\bk )\right]_{\bk\in I_K^d}$.}

  \caption{Inversion of $\bm D^{\mathrm{tot}}$} 
  \label{alg:bw}
\end{algorithm}

The DFT \eqref{eq:dft} and its inverse \eqref{eq:idft} can be computed
efficiently with the well-known \emph{fast Fourier transform} (FFT) in
$\mathcal O(N^{d-1} \log N)$ operations.  Furthermore, there are fast
algorithms for the computation of the NDFT \eqref{eq:ndft} known as
\emph{nonuniform fast Fourier transform} (NFFT), which provides an
arbitrarily tight approximation of the NDFT and uses
$\mathcal{O}(K^d\log K)$ instead of the $\mathcal O(K^{2d})$
operations of the naive matrix-vector multiplication
\cite{bey95,duro93,st97}.  For our numerical simulations, we use the
implementation of the NFFT \cite{KeKuPo09} provided in the software
package \cite{nfft3}, which also includes a fast algorithm for the
adjoint NDFT with the same arithmetic complexity.

The NDFT \eqref{eq:ndft} may also be interpreted
as the (discretized) Fourier integral operator
$
g(\bl, m)
=
\sum_{\bk\in\mathcal I_K^d} \e^{\i\, \Phi((\bl,m),\bk)} f_\bk
$
with $\Phi((\bl,m),\bk) = -2\ls K^{-1}\, (R_{t_m}\bh(\by'_\bl)) \cdot \bk$,
such operators are more generally studied in \cite{CanDemYin09}.
This allows to apply fast butterfly algorithms \cite{LiYaYi15}, which are of the same asymptotic complexity as the NFFT.

\section{Inverse NDFT as inverse problem}
\label{sec:indft}

The approximate inversion of the NDFT \eqref{eq:ndft} is the most
crucial part for the reconstruction in \autoref{alg:bw}.  Let
$\bg\in\R^{MN^{d-1}}$ be given.  An (approximate) solution
$\bm f\in\R^{K^d}$ of
\begin{equation}\label{eq:indft}
  \ndft \bm f = \bm g
\end{equation} 
is called an \emph{inverse NDFT} of $\bm g$.
There is no exact inversion formula known in general.
Depending on the distribution of the nodes $\mathcal Y_{M,N}$, the inverse NDFT \eqref{eq:indft} might not have a unique solution or be ill-conditioned \cite[Sect.\ 7.6.2]{PlPoStTa18}.

\subsection{Discrete backpropagation}
\label{sec:backprop}

The idea behind the backpropagation \cite{Dev82} is to apply a
quadrature rule to the (continuous) inverse Fourier transform
restricted to the set $\mathcal Y$ in \eqref{eq:Y}, where values
are available by the Fourier diffraction theorem \eqref{eq:recon}.
The \emph{continuous backpropagation} is defined by
\begin{equation} \label{eq:bp}
(2\pi)^{-d/2}
\int_{\mathcal Y} \mathcal F[f](\by)\, \e^{\i\, \by\cdot\bx}
\dd \by
,\qquad \bx\in\R^d,
\end{equation}
which is simply the inverse Fourier transform applied to
$\bm1_{\mathcal Y}\,\mathcal F[f] $, where
$\bm1_{\mathcal Y}(\by)\coloneqq1$ if $\by\in\mathcal Y$ and zero
elsewhere.  Applying the substitution
$ (t,\bz') \mapsto \by = R_t \bh(\bz') $ and a quadrature rule with
nodes $(t_m,\bz_\bn')$ and weights
$\bm w = [w_{m,\bn}]_{m\in\mathcal J_M, \bn\in\mathcal I_N^{d-1}} \in
\R_{>0}^{MN^{d-1}}$ to \eqref{eq:bp}, we arrive at the \emph{discrete
  backpropagation (BP)} $\bm B$ of $\bm g$, namely
\begin{equation} \label{eq:bp-d}
[\bm B\bg]_\bk
\coloneqq
\sum_{m\in\mathcal J_M} 
\sum_{\bn\in\mathcal I_N^{d-1}}
w_{m,\bn}\,
g_{m,\bn}\,
\e^{-\i\, (\bk\frac{2\ls}{K})\cdot (R_{t_m} \bh(\bz'_\bn))}
,\qquad\bk\in \mathcal I_K^d.
\end{equation}
For many situations in ODT, the weights $\bw$ are known analytically
\cite[Sect.\ 4]{KirQueRitSchSet21}.  For instance, if the 3D object
makes a full turn around the $x_1$ axis with rotation angle
$t\in[0,2\pi]$, we have the nodes \eqref{eq:z'} and weights
\begin{equation} \label{eq:w-e1}
w_{m,\bn} 
= \frac{k_0 \abs{[\bz'_\bn]_2}}{2\,\kappa(\bz'_\bn)}
\, \frac{2\pi^2 k_0^2}{MN^{2}}
,
\qquad m\in\mathcal J_M,\ \bn\in\mathcal I_N^2.
\end{equation}
The discrete backpropagation is an adjoint NDFT, i.e.\
$ \bm B\bg = \ndft^* (\bw \odot \bg), $ which can be efficiently
evaluated via the adjoint NFFT.

\subsection{Conjugate gradient method}

\label{sec:cg}

\begin{algorithm}[b!]
  \KwIn{$\bg = [g_{m,\bn}]_{m\in\mathcal J_M, \bn\in\mathcal I_N} \in
    \C^{MN^{d-1}}$, $\bm w \in \R_{>0}^{MN^{d-1}}$,
    $\bm f^{(0)}\in\R^{K^d}$, $J_{\mathrm{CG}}\in\N$.  }

  $\bm r^{(0)} \coloneqq \bg$\;
  $\bm p^{(0)} \coloneqq \bm s^{(0)} \coloneqq \re[\ndft^*(w\odot \br^{(0)})] $\;
  
  \For{$j = 0,1,\dots,J_{\mathrm{CG}}$}{
    $\bm q^{(j)} \coloneqq \ndft \bm p^{(j)} $\;
    $\alpha_j \coloneqq (\bm s^{(j)})^*\,\bm s^{(j)} / (\bm q^{(j)})^*\, (\bm w\odot\bm q^{(j)})$\;
    $\bm f^{(j+1)} \coloneqq \bm f^{(j)} + \alpha_j\, \bm p^{(j)} $\;
    $\br^{(j+1)} \coloneqq \br^{(j)}-\alpha_j\, \bm q^{(j)}$\;
    $\bm s^{(j+1)} \coloneqq \re [\ndft^* (\bm w\odot\bm r^{(j+1)})]$\;
    $\beta_j \coloneqq (\bm s^{(j+1)})^*\,\bm s^{(j+1)} / (\bm s^{(j)})^*\,\bm s^{(j)}$\;
    $\bm p^{(j+1)} \coloneqq \bm s^{(j+1)} + \beta_j\,\bm p^{(j+1)}$\;
  }
  
  \KwOut{Approximate solution $\bm f^{(J_{\mathrm{CG}})} \in \R^{K^d}$.}
  
  \caption{CG inversion of the NDFT} 
  \label{alg:cg}
\end{algorithm}

Equation \eqref{eq:indft} is a linear equation system that may be
solved in the least-square sense by applying the iterative
\emph{conjugate gradient} (CG) method.  This approach has been
successfully applied in ODT \cite{FauKirQueSchSet21,KirQueRitSchSet21}
as well as in X-ray imaging \cite{PoSt01IMA}, magnetic resonance
imaging \cite{KnKuPo}, or spherical tomography \cite{HiQu15}.  We
consider the overdetermined case $MN^{d-1} \ge K^d$.  The
underdetermined case is similar and covered in \cite{kupo04}.  The
idea is to find a solution $\bm f\in\R^{K^d}$ of
\begin{equation} \label{eq:indft-min}
  \minimize_{\bm f \in \mathbb R^{K^d}} \qquad
  \norm{\ndft\bm f - \bm g}_{2,\zb w}^2
\end{equation}
with some weights
$\bm w = [w_{m,\bn}]_{m\in\mathcal J_M, \bn\in\mathcal I_N^{d-1}} \in
\R_{>0}^{MN^{d-1}}$.  The classical CG method corresponds to the
constant weights $w_{m,\bn}\equiv1$.  Obviously, any exact solution
$\bm f$ of \eqref{eq:indft} is a minimizer of \eqref{eq:indft-min}.
We apply the CG method on the normal equation
$$
\re[\ndft^* (\bm w \odot \ndft\bm f)]
=
\re[\ndft^* (\bm w\odot \bm g)],
$$
where $\bm f$ is assumed to be real-valued.  This approach results in
\autoref{alg:cg}, which is known as CGNE (CG on the normal equations),
CGNR (CG minimizing the norm of the residual) or CGLS (CG minimizing
least squares).  The iterative method requires the evaluation of
$\ndft$ and $\ndft^*$ in each iteration.  For a more detailed
description of the algorithm, we refer to \cite{KnKuPo,Hankebook}.

Good weights $\bm w$ ensure that $\ndft^* (\bm w \odot \bm g)$ is
close to the desired solution $\bm f$.  In this case, the 
CG usually converges faster.  One possibility to choose
$\bm w$ are the quadrature weights \eqref{eq:w-e1} with respect to
$\mathcal Y_{M,N}$ \cite{KnKuPo}.  The CG algorithm is a simple
regularization method, where the number of iterations acts as
regularization parameter \cite[Sect.\ 7]{EnHaNe96}.

\subsection{Total variation regularization}

\label{sec:tv}

In order to stabilize the inversion of the NDFT, we regularize
\eqref{eq:indft} by the total variation (TV) semi-norm.  For a
discrete image $\bm f \in \mathbb R^{K^d}$, the \emph{total variation}
becomes
\begin{equation*}
  \mathrm{TV}(\bm f)
  \coloneqq \sum_{\bm k \in \mathcal I_K^d}
  \| ( \grad \bm f )_{\bm k} \|_2
  = \| \grad \bm f \|_{1,2},
\end{equation*}
see for instance Bredies \& Lorenz \cite{BreLor18}, where the discrete
gradient $\grad \colon \mathbb R^{K^d} \to \mathbb R^{K^d \times d}$
yields a discrete vector field by applying the finite forward
differences with respect to each axis.  More precisely, the discrete
gradient is defined by
\begin{equation} \label{eq:grad}
  (\grad \bm f)_{\bm k,\ell} =
  (\partial_\ell \bm f)_{\bm k}
  \qquad\text{with}\qquad
  (\partial_\ell \bm f)_{\bm k}
  \coloneqq
  \begin{cases}
    f_{\bm k + \bm e_\ell} - f_{\bm k},
    & k_\ell < \frac K2 - 1,
    \\
    0,
    &  k_\ell = \frac K2 - 1,
  \end{cases}
\end{equation}
for $\bm k \in \mathcal I_K^d$ and $\ell \in \mathcal J_d$, where
$\bm e_\ell$ corresponds to the $\ell$th unit vector. The adjoint of
the discrete gradient is the negative discrete divergence, i.e.\
$\grad^* = - \Div$.  Using finite backward differences, we define the
discrete divergence by
\begin{equation*}
  \Div \bm y
  \coloneqq \sum_{\ell=1}^d \tilde\partial_\ell \bm y_{\cdot,\ell}
  \qquad\text{with}\qquad
  (\tilde\partial_\ell \bm f)_{\bm k}
  \coloneqq
  \begin{cases}
    f_{\bm k},
    & k_\ell = -\tfrac K2,
    \\
    f_{\bm k} - f_{\bm k - \bm e_\ell},
    & -\tfrac K2 < k_\ell < \tfrac K2 - 1,
    \\
    -f_{\bm k-\be_\ell},
    & k_\ell = \tfrac K2 - 1,
  \end{cases}
\end{equation*}
where $\bm y_{\cdot, \ell}$ denotes the $\ell$th components of the
discrete vector field $\bm y \in \mathbb R^{K^d \times d}$.

The TV regularization is well known to promote
cartoon-like images such as the objects we expect in our specific
application.  Besides this a priori information, we know that the
images representing the physical object are real-valued and
non-negative. Incorporating this knowledge, we propose to compute the
NDFT by minimizing the variational formulation
\begin{equation}
  \label{eq:var-prob}
  \minimize_{\bm f \in \mathbb R^{K^d}} \qquad
  \underbracket{\chi_{\mathbb
      R^{K^d}_{\ge0}}(\bm f)}_{\eqqcolon F(\bm f)}
  + \underbracket{\tfrac12 \|\bm F_{\mathrm{NDFT}}
    (\bm f) -
    \bm g\|^2_{2,\bm w}}_{\eqqcolon G(\bm f)} +
  \underbracket{\lambda \| \grad \bm f \|_{1,2}}_{\eqqcolon H (L\bm f)},
\end{equation}
where $L\coloneqq\grad$ and
$\chi_{\mathbb R_{\ge 0}^{K^d}}$ is the characteristic function over
the non-negative orthant given by
\begin{equation*}
  \chi_{\mathbb R^{K^d}_{\ge 0}}(\bm f)
  \coloneqq
  \begin{cases}
    0, & \text{if} \; f_{\bm k} \ge 0 \; \text{for all} \; k \in
    \mathcal I_K^d, \\
    +\infty, & \text{otherwise}.
  \end{cases}
\end{equation*}

A convex minimization problem like \eqref{eq:var-prob} may be
numerically solved with the primal-dual (PD) iteration by Chambolle \& Pock
\cite{ChaCasCreNovPoc10}, which consists in
\begin{equation*}
  \begin{aligned}
    \bm f^{(j+1)}
    &\coloneqq \prox_{\tau F}
    \bigl(\bm f^{(j)} - \tau (\nabla G(\bm f^{(j)}) + L^* \bm y^{(j)})\bigr),
    \\
    \bm y^{(j+1)}
    &\coloneqq \prox_{\sigma H^*}
    \bigl(\bm y^{(j)} + \sigma L (2 \bm f^{(j+1)} - \bm f^{(j)})\bigr)
  \end{aligned}
\end{equation*}
for appropriate $\tau, \sigma > 0$.  Here, $\nabla G$ denotes the
gradient of the differentiable function $G$ in contrast to the
discrete gradient \eqref{eq:grad} of a multidimensional signal.  More
generally, the \emph{proximal mapping} of a convex function
$h \colon \mathbb R^K \to \mathbb R \cup \{-\infty,+\infty\}$ is
defined by
\begin{equation}
  \label{eq:prox}
  \prox_h(\bm x)
  \coloneqq
  \argmin_{\bm y \in \mathbb R^K} h(\bm y) + \tfrac12 \, \|\bm x  - \bm y \|_2^2.
\end{equation}
If $h$ is additionally proper, i.e.\ not everywhere $+\infty$ and
never $-\infty$, and lower semi-continuous, then the minimizer in
\eqref{eq:prox} is unique, and the proximal mapping becomes single
valued.  Further, the \emph{Fenchel} or \emph{convex conjugate} of the
convex function $h$ is defined by
\begin{equation*}
  h^*(\bm x) \coloneqq
  \sup_{\bm y \in \mathbb R^K} \langle \bm x , \bm y \rangle - h(\bm y).
\end{equation*}

Since the proximal mapping of an indicator function is the projection to
its support \cite[Thm.~3.14]{BauCom11}, i.e.\
\begin{equation*}
  \prox_{\mathbb R_{\ge0}^{K^d}} (\bm x)
  = \proj_{\mathbb R_{\ge0}^{K^d}} (\bm x)
  \coloneqq \argmin_{\bm y \in \mathbb R_{\ge 0}^{K^d}} \| \bm x - \bm
  y \|_2^2,
\end{equation*}
for our specific formulation \eqref{eq:var-prob}, the PD
iteration becomes
\begin{equation}
  \label{eq:pri-dual}
  \begin{aligned}
    \bm f^{(j+1)}
    &\coloneqq \proj_{\mathbb R^{K^d}_{\ge0}}
    (\bm f^{(j)} - \tau (\re[\bm F_{\mathrm{NDFT}}^*(\bm w \odot (\bm
    F_{\mathrm{NDFT}} 
    \bm f^{(j)} -
    \bm g))] - \Div \bm y^{(j)})),
    \\
    \bm y^{(j+1)}
    &\coloneqq \prox_{\sigma \lambda \|\cdot\|^*_{1,2}}
    (\bm y^{(j)} + \sigma \grad (2 \bm f^{(j+1)} - \bm f^{(j)})).
  \end{aligned}
\end{equation}
The remaining proximal mapping of the discrete TV semi-norm may be
computed by exploiting
\begin{equation*}
  \prox_{\sigma \lambda \|\cdot\|^*_{1,2}}(\bm y)
  = \bm y - \sigma \prox_{\frac\lambda\sigma
    \|\cdot\|_{1,2}}(\tfrac1\sigma \bm y),
\end{equation*}
cf.\ \cite[Thm.~14.3]{BauCom11}, and
\begin{equation*}
  (\prox_{\rho \| \cdot \|_{1,2}}(\bm y))_{\bm k, \ell}
  = \bigl( 1 - \tfrac\rho{\max\{\| \bm y_{\bm k, \cdot} \|_2, \rho\}}
  \bigr) \, y_{\bm k, \ell}
\end{equation*}
for some $\rho>0$,
see for example \cite{ChiChoComPes} and references therein.

The central pitfall of applying the PD iteration in practice
is the selection of appropriate primal and dual step sizes $\tau$ and
$\sigma$.  In the literature, there exist several approaches to choose
the step sizes adaptively \cite{YokHon17,GolLiYua15}.  In our
numerical experiments, we rely on the step size rules proposed in
\cite{YokHon17}, which have been successfully applied to TV denoising
in image processing.  The main idea is to compute the primal and dual
residuals
\begin{equation}
  \label{eq:resi}
  \begin{aligned}
    \bm p^{(j+1)}
    &\coloneqq \tfrac1\tau \, (\bm f^{(j)} -
    \bm f^{(j+1)}) - \re[\bm F_{\mathrm{NDFT}}^*(\bm w \odot \bm F_{\mathrm{NDFT}}
    (\bm f^{(j)} - \bm f^{(j+1)}))] + \Div(\bm y^{(j)} -
    \bm y^{(j+1)}),
    \\
    \bm d^{(j+1)}
    &\coloneqq \tfrac1\sigma \, (\bm y^{(j)} - \bm y^{(j+1)})
    - \grad (\bm f^{(j)} - \bm f^{(j+1)}),
  \end{aligned}
\end{equation}
where we added the residual with respect to the gradient $\nabla G$.
Further, we require the values
\begin{equation}
  \label{eq:b-con}
  \begin{aligned}
    \alpha^{(j+1)}
    &\coloneqq \tfrac{\| \bm f^{(j+1)}
      \|_2}{\|\bm y^{(j+1)}\|_{2,2}},
    \\
    \omega_{\mathrm P}^{(j+1)}
    &\coloneqq \frac{\langle \bm f^{(j)}
      - \bm f^{(j+1)}, \bm p^{(j+1)} \rangle}{\|\bm
      f^{(j)} - \bm f^{(j+1)}\|_2 \|\bm p^{(j+1)} \|_2},
    \\
    \omega_{\mathrm D}^{(j+1)}
    &\coloneqq \frac{\langle \bm y^{(j)}
      - \bm y^{(j+1)}, \bm d^{(j+1)} \rangle}{\|\bm
      y^{(j)} - \bm y^{(j+1)}\|_{2,2} \|\bm d^{(j+1)} \|_{2,2}},
  \end{aligned}
\end{equation}
where the occurring inner products are defined via
\begin{equation*}
  \langle \bm f, \bm g \rangle
  \coloneqq \sum_{\bm k \in \mathcal I_K^d} f_{\bm k} g_{\bm
    k}
  \qquad\text{and}\qquad
  \langle \bm y, \bm z \rangle
  \coloneqq \sum_{\bm k \in \mathcal I_K^d} \langle \bm
  y_{\bm k, \cdot}, \bm z_{\bm k, \cdot} \rangle.
\end{equation*}

\begin{algorithm}[bt]
  \KwIn{$\bm g = \bm F_{\mathrm{NDFT}} \bm f$,
    $\bm w \in \R_{>0}^{MN^{d-1}}$,
    $\lambda>0$,
    $\bm f^{(0)} \in \mathbb R^{K^d}_{\ge0}$,
    $\bm y^{(0)} \in \mathbb R^{K^d \times d}$,
    $\tau^{(0)}>0$, $\sigma^{(0)}>0$, $\rho = 0.005$, $c = 0.9$,
    $\beta = 1.5$, $\zeta = 0.25$, $J_{\mathrm{PD}}\in\N$.}
  
  \For{$j=0,1,2,\dots,J_{\mathrm{PD}}$}{
    Compute $\bm f^{(j+1)}$, $\bm y^{(j+1)}$ by \eqref{eq:pri-dual}\;
    Compute $\bm p^{(j+1)}$, $\bm d^{(j+1)}$ by \eqref{eq:resi}\;
    Compute $\alpha^{(j+1)}$, $\omega_{\mathrm P}^{(j+1)}$,
    $\omega_{\mathrm D}^{(j+1)}$ by \eqref{eq:b-con}\;
    \lIf{$\omega_{\mathrm P}^{(j+1)} > c$}{$\tau^{(j)} \coloneqq
      \beta \tau^{(j)}$}
    \lIf{$\omega_{\mathrm P}^{(j+1)} < 0$}{$\tau^{(j)} \coloneqq
      \zeta \tau^{(j)}$}
    \lIf{$\omega_{\mathrm D}^{(j+1)} > c$}{$\sigma^{(j)} \coloneqq
      \beta \sigma^{(j)}$}
    \lIf{$\omega_{\mathrm D}^{(j+1)} < 0$}{$\sigma^{(j)} \coloneqq
      \zeta \sigma^{(j)}$}
    Balancing: $\tau^{(j+1)} \coloneqq (\alpha^{(j+1)})^\rho \tau^{(j)}$,
    $\sigma^{(j+1)} \coloneqq \sigma^{(j)} / (\alpha^{(j+1)})^{\rho}$\;
  }
  \KwOut{Approximate scattering potential $\bm f^{(J_{\mathrm{PD}})}$.}
  \caption{PD inversion of the NDFT with TV} 
  \label{alg:pd-back}
\end{algorithm}

Now, if $\bm f^{(j)} - \bm f^{(j+1)}$ is aligned with primal residual
$\bm p^{(j+1)}$, we increase the primal step size $\tau$; and if
$\bm f^{(j)} - \bm f^{(j+1)}$ is opposed to $\bm p^{(j+1)}$, we
decrease $\tau$.  We proceed analogously for the dual step size
$\sigma$.  On the basis of the observation that a large $\tau$ and a
small $\sigma$ decrease the primal residual at the cost of the dual
residual \cite{GolLiYua15}, Yokata \& Hontani \cite{YokHon17} propose
to balance both step sizes by multiplying $\tau$ with a power of
$\alpha^{(j+1)}$ and dividing $\sigma$ by the same power of
$\alpha^{(j+1)}$, which is greater than one if the primal residual
dominates; or less than one if the dual residual dominates.
Summarized, we apply the PD iteration with the parameter updates in
\autoref{alg:pd-back}, which we call the PD method (with TV) for ODT.
A more detailed motivation behind the acceleration, backtracking and
balancing of the step sizes may be found in \cite{YokHon17}.

\subsection{Numerical backward transform}
\label{sec:numer-backw-transf}

We numerically compare the three reconstruction approaches based on
\autoref{alg:bw}, where the inverse NDFT step is realized via
\begin{enumerate}[i)]
\item the BP method \eqref{eq:bp-d},
\item the CG method in \autoref{alg:cg}, and
\item the PD method with TV in \autoref{alg:pd-back}.
\end{enumerate}
Additionally, 
we apply TV denoising to the reconstructed images of i) and ii).
The \emph{TV denoising} (TVd) corresponds to
finding a minimizer of 
$$
\minimize_{\tilde{\bm f} \in \mathbb R^{K^d}}
\qquad
\chi_{\R_{\ge0}^{K^d}}(\tilde{\zb f})
+\tfrac12\norm{\smash{\tilde{\zb f} - \zb f}}_2^2
+\lambda \norm{\smash{\grad \tilde{\zb f}}}_{1,2},
$$
i.e., \eqref{eq:var-prob} with the identity instead of
$\bm F_{\mathrm{NDFT}}$.  In order to avoid the so-called ``inverse
crime,'' we compute the data via a direct discretization of the
convolution \eqref{eq:conv}. 
For some of the tests, we also disturb the
measurements by additive noise; so in this case we have only access to
$\bm D^{\mathrm{tot}} \bm f + \zb\varepsilon$ with
$\zb\varepsilon\in\C^{M N^{d-1}}$ corresponding to the noise level
$\norm{\zb\varepsilon}_2/\norm{\bm D^{\mathrm{tot}} \bm f}_2$.

All numerical tests were performed on a
standard PC with an 8-core Intel Core i7-10700 and 32\,GB memory using
Matlab with the NFFT library \cite{nfft3}.
Our code is available in the \texttt{FourierODT} toolbox.\footnote{\url{https://github.com/michaelquellmalz/FourierODT}}
We compare the reconstruction quality by means of the \emph{structural similarity index measure} (SSIM) \cite{ssim}
and the
\emph{peak signal-to-noise ratio}
$$
\operatorname{PSNR}(\zb f,\zb g)
\coloneqq 10 \log_{10} \frac{\max_{\br\in\mathcal I_K^d} \abs{f_\bk}^2}{K^{-d} \sum_{\bk\in\mathcal I_K^d} \abs{f_\bk - g_\bk}^2},
\qquad \zb f, \zb g \in \R^{K^d}.
$$

\paragraph{Two-dimensional function.}

For a 2D simulation, we employ the discretization parameters
$K=N=M=240$ as well as
$\rM=40$ and $\lM=60$.  
We choose $\ls = K/(4\sqrt2)\approx42$,
which ensures that $\mathcal{Y_{M,N}}$ is inside one peridicity interval of the NDFT.
We fix the wave number $k_0=2\pi$,
which means that all spatial measurements are in multiples of the
wavelength of the incident field $\ui$.  The object rotates a full
turn around the center.  The object $\bm f$ consists of different
convex and concave shapes, see \autoref{fig:heart1-f}. 
The discrete sinogram $\bm D^{\mathrm{tot}} \bm f$ and the
Fourier space data approximating $\mathcal F[f]$ on the set
\eqref{eq:Y_MN}, which is contained in a disk of radius $\sqrt2\,k_0$,
are shown in \autoref{fig:heart1-sinogram}.  

\begin{table}[htp]\centering
  \begin{tabular}{c c c c c c}
    & BP & BP\,\&\,TVd & CG & CG\,\&\,TVd & PD
    \\\hline
    PSNR
    & 31.22 & 36.17 & 39.61 & 40.12 & 41.59
    \\\hline SSIM
    & 0.388 & 0.991 & 0.983 & 0.990 & 0.988
    \\\hline Time (sec)
    & 0.01 & 0.12 & 0.18 & 0.12 & 0.28
  \end{tabular}
  \caption{Error for the reconstruction of the function of \autoref{fig:heart1-f} with exact data.
  \label{tab:heart1}}
\end{table}

\tikzset{font=\tiny}
\newcommand{\modelwidth} {3.6cm}
\newcommand{\smallmodelwidth} {2.4cm}
\pgfplotsset{
  colormap={parula}{
    rgb255=(53,42,135)
    rgb255=(15,92,221)
    rgb255=(18,125,216)
    rgb255=(7,156,207)
    rgb255=(21,177,180)
    rgb255=(89,189,140)
    rgb255=(165,190,107)
    rgb255=(225,185,82)
    rgb255=(252,206,46)
    rgb255=(249,251,14)}}

\pgfmathsetmacro{\xmin}{-42.426407}
\pgfmathsetmacro{\xmax}{42.072853}

\begin{figure}[htp]\centering
  \begin{subfigure}[t]{.33\textwidth}\centering
    \begin{tikzpicture}
      \begin{axis}[
        width=\modelwidth,
        height=\modelwidth,
        enlargelimits=false,
        scale only axis,
        axis on top,
        colorbar,colorbar style={
          width=.15cm, xshift=-0.6em, 
          point meta min=0.8914,point meta max=1.9101,
        },
        ]
        \addplot graphics [
        xmin=0, xmax=6.28,
        ymin=-60, ymax=60,
        ] {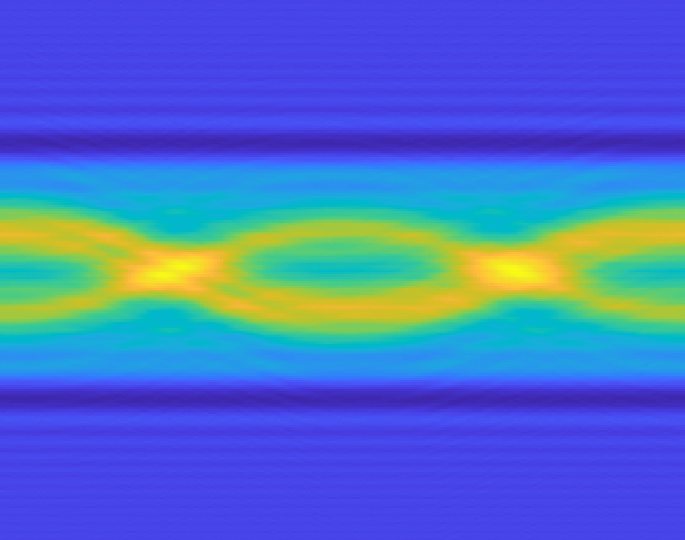};
      \end{axis}
    \end{tikzpicture}
    \caption{Sinogram $\abs{\utot_t(x)}$ for the function of \autoref{fig:heart1-f}}
  \end{subfigure}\hfill
  \begin{subfigure}[t]{.33\textwidth}\centering
    \begin{tikzpicture}
      \begin{axis}[
        width=\modelwidth,
        height=\modelwidth,
        enlargelimits=false,
        scale only axis,
        axis on top,
        colorbar,colorbar style={
          width=.15cm, xshift=-0.6em, 
          point meta min=0.7475, point meta max=2.0076,
        },
        ]
        \addplot graphics [
        xmin=0, xmax=6.28,
        ymin=-60, ymax=60,
        ] {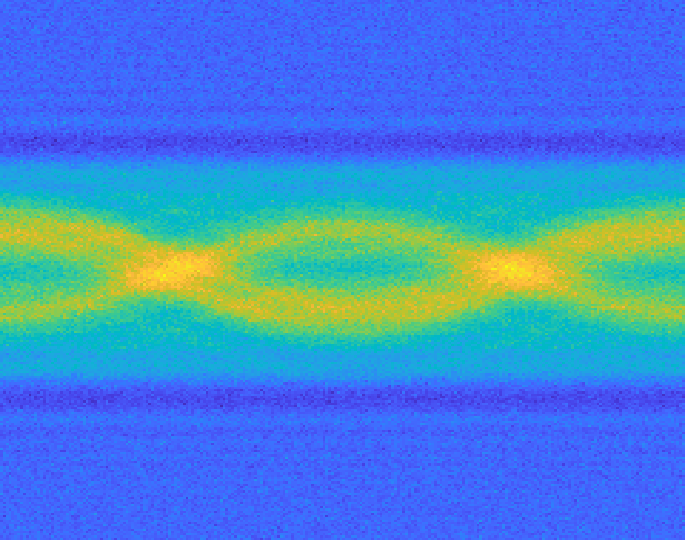};
      \end{axis}
    \end{tikzpicture}
    \caption{Sinogram with 5\,\% Gaussian white noise}
  \end{subfigure}\hfill
  \begin{subfigure}[t]{.33\textwidth}\centering
    \begin{tikzpicture}
      \begin{axis}[width=\modelwidth, height=\modelwidth, enlarge y limits=0.1, axis x line*=bottom, axis y line*=left, legend pos=outer north east, legend style={cells={anchor=west}}, scale only axis, 
        ]
        \addplot+[only marks,scatter,mark size=0.5,point meta={\thisrowno{2}}, scatter/use mapped color=
        {draw opacity=0,fill=mapped color}
        ] table[x index=0,y index=1, meta index=2] {heart1.txt};
      \end{axis}
    \end{tikzpicture}
    \caption{Fourier space data on $\mathcal Y_{M,N}$ for a full turn
    }
  \end{subfigure}
  \caption{Visualization of the data and the respective Fourier space coverage.
  \label{fig:heart1-sinogram}}
\end{figure}
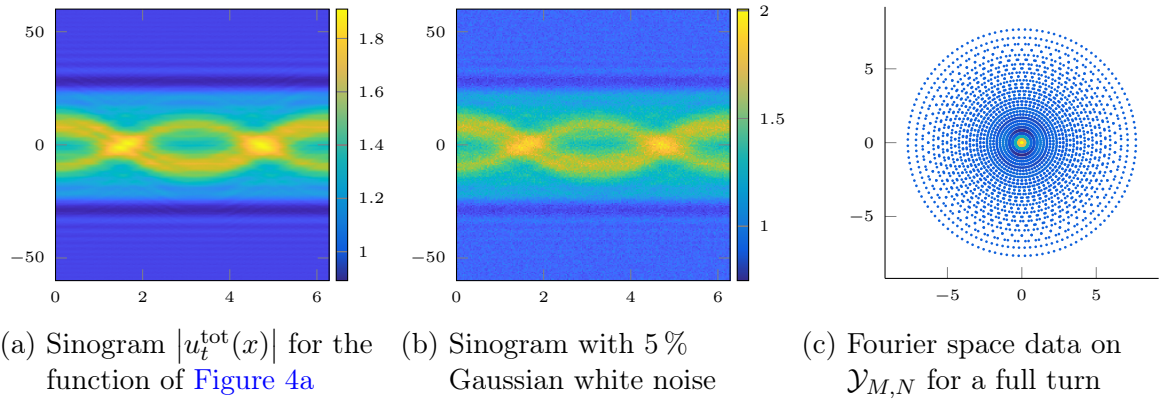

\begin{figure}[htp]\centering
\pgfmathsetmacro{\cmin} {-0.05} \pgfmathsetmacro{\cmax} {0.55}
\begin{subfigure}[t]{.33\textwidth}
  \begin{tikzpicture}
    \begin{axis}[
      width=\modelwidth,
      height=\modelwidth,
      enlargelimits=false,
      scale only axis,
      axis on top,
      colorbar,colorbar style={
        width=.15cm, xshift=-0.9em, 
        point meta min=\cmin,point meta max=\cmax,
      },
      ]
      \addplot graphics [
      xmin=\xmin, xmax=\xmax,  ymin=\xmin, ymax=\xmax,
      ] {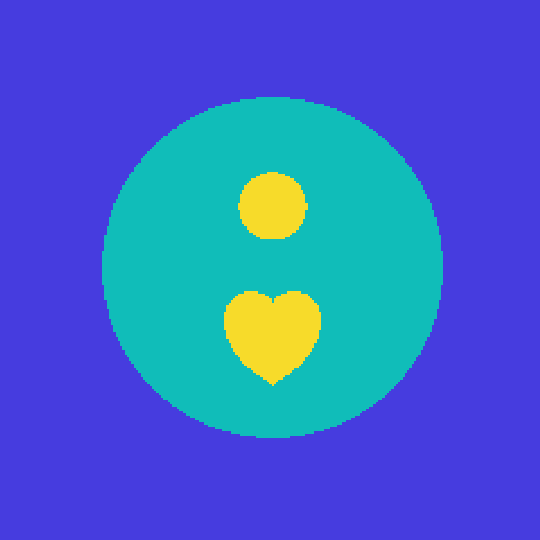};
    \end{axis}
  \end{tikzpicture}
  \caption{Ground truth $\bm f$ \label{fig:heart1-f}}
\end{subfigure}\hfill
\begin{subfigure}[t]{.33\textwidth}\centering
  \begin{tikzpicture}
    \begin{axis}[
      width=\modelwidth,
      height=\modelwidth,
      enlargelimits=false,
      scale only axis,
      axis on top,
      colorbar,colorbar style={
        width=.15cm, xshift=-0.9em, 
        point meta min=\cmin,point meta max=\cmax,
      },
      ]
      \addplot graphics [
      xmin=\xmin, xmax=\xmax,  ymin=\xmin, ymax=\xmax,
      ] {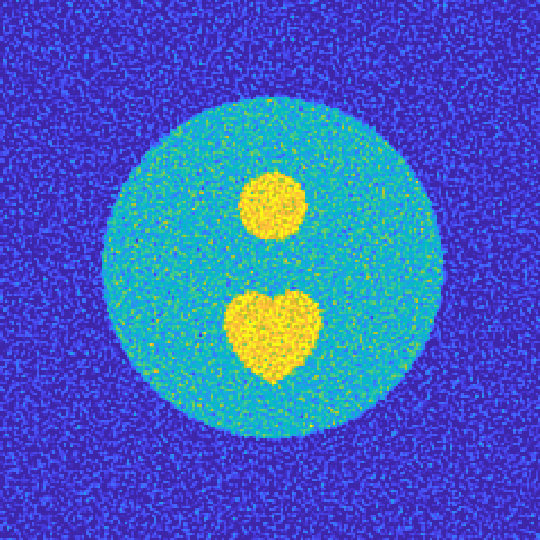};
    \end{axis}
  \end{tikzpicture}
  \caption{BP reconstruction\\PSNR 23.38, SSIM 0.125}
\end{subfigure}\hfill
\begin{subfigure}[t]{.33\textwidth}\centering
  \begin{tikzpicture}
    \begin{axis}[
      width=\modelwidth,
      height=\modelwidth,
      enlargelimits=false,
      scale only axis,
      axis on top,
      colorbar,colorbar style={
        width=.15cm, xshift=-0.9em, 
        point meta min=\cmin,point meta max=\cmax,
      },
      ]
      \addplot graphics [
      xmin=\xmin, xmax=\xmax,  ymin=\xmin, ymax=\xmax,
      ] {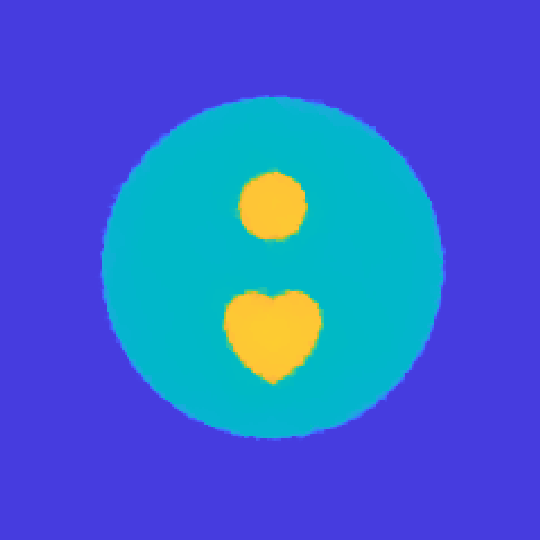};
    \end{axis}
  \end{tikzpicture}
  \caption{BP and TVd ($\lambda=0.1$)\\PSNR 34.71, SSIM 0.985}
\end{subfigure}

\begin{subfigure}[t]{.33\textwidth}\centering
  \begin{tikzpicture}
    \begin{axis}[
      width=\modelwidth,
      height=\modelwidth,
      enlargelimits=false,
      scale only axis,
      axis on top,
      colorbar,colorbar style={
        width=.15cm, xshift=-0.9em, 
        point meta min=\cmin,point meta max=\cmax,
      },
      ]
      \addplot graphics [
      xmin=\xmin, xmax=\xmax,  ymin=\xmin, ymax=\xmax,
      ] {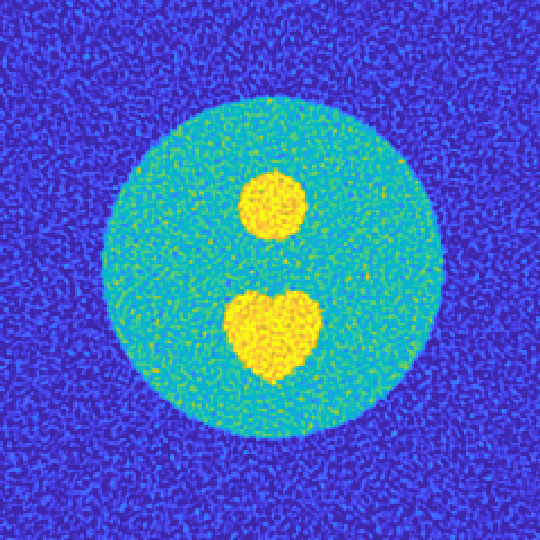};
    \end{axis}
  \end{tikzpicture}
  \caption{CG reconstruction\\PSNR 26.03, SSIM 0.234}
\end{subfigure}\hfill
\begin{subfigure}[t]{.33\textwidth}\centering
  \begin{tikzpicture}
    \begin{axis}[
      width=\modelwidth,
      height=\modelwidth,
      enlargelimits=false,
      scale only axis,
      axis on top,
      colorbar,colorbar style={
        width=.15cm, xshift=-0.9em, 
        point meta min=\cmin,point meta max=\cmax,
      },
      ]
      \addplot graphics [
      xmin=\xmin, xmax=\xmax,  ymin=\xmin, ymax=\xmax,
      ] {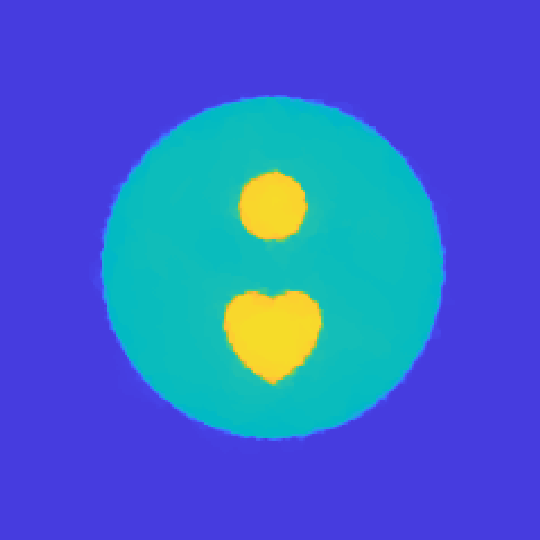};
    \end{axis}
  \end{tikzpicture}
  \caption{CG and TVd ($\lambda=0.1$)\\ PSNR 38.05, SSIM 0.983}
\end{subfigure}\hfill
\begin{subfigure}[t]{.33\textwidth}\centering
  \begin{tikzpicture}
    \begin{axis}[
      width=\modelwidth,
      height=\modelwidth,
      enlargelimits=false,
      scale only axis,
      axis on top,
      colorbar,colorbar style={
        width=.15cm, xshift=-0.9em, 
        point meta min=\cmin,point meta max=\cmax,
      },
      ]
      \addplot graphics [
      xmin=\xmin, xmax=\xmax,  ymin=\xmin, ymax=\xmax,
      ] {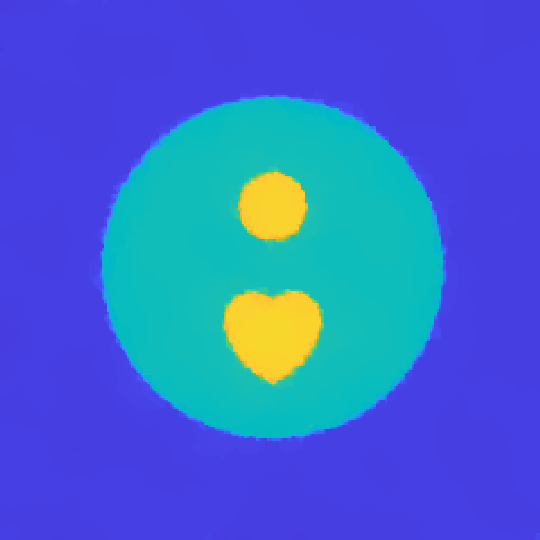};
    \end{axis}
  \end{tikzpicture}
  \caption{PD with TV ($\lambda=0.1$)\\PSNR 37.62, SSIM 0.872 \label{fig:heart1-noise-tv}}
\end{subfigure}%
\caption{Inversion of $\bm D^{\mathrm{tot}} \bm f$ with 5\% Gaussian
  noise using BP in \eqref{eq:bp-d}, CG in Alg.~\ref{alg:cg}, and PD
  with TV in Alg.~\ref{alg:pd-back}.  The BP and CG results are
  TV-denoised.
  \label{fig:heart1-noise}}
\end{figure}

For the reconstructions, we use $J_{\mathrm{CG}}=20$
iterations for CG in \autoref{alg:cg} and $J_{\mathrm{PD}}=50$
iterations for PD with TV in \autoref{alg:pd-back}.  Note that the
BP and the CG reconstruction only assume that $\bm f$ is
real-valued, whereas PD with TV incorporates the non-negativity of
$\bm f$.
The reconstruction error for exact data and the computation times are shown in \autoref{tab:heart1}
(all reconstructions look very good and are visually almost indistinguishable).
\autoref{fig:heart1-noise} shows the reconstructions for 5\% complex
Gaussian noise.  
Different from above, the CG method is performed with $J_{\mathrm{CG}}=5$
iterations, which acts as regularization against the noise. The
advantage of the TV regularization is clearly visible.

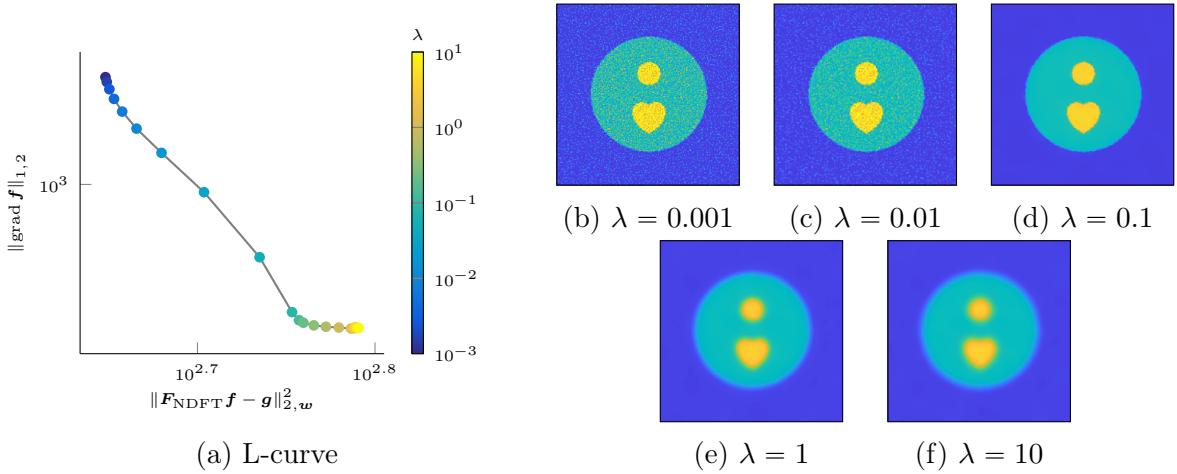
\begin{figure}[tp]\centering
  \pgfmathsetmacro{\cmin} {-0.05} \pgfmathsetmacro{\cmax} {0.55}
\begin{subfigure}[b]{.45\textwidth}
  \begin{tikzpicture}
    \pgfmathdeclarefunction{lg10}{1}{\pgfmathparse{ln(#1)/ln(10)}}
    \begin{loglogaxis}[width=4cm, height=4cm, enlarge y limits=0.1, axis x line*=bottom, axis y line*=left, legend pos=outer north east, legend style={cells={anchor=west}}, scale only axis, 
      colorbar, colorbar style={width=.15cm, xshift=-0.6em,
      ytick={-3,-2,-1,0,1},yticklabels={$10^{-3}$,$10^{-2}$,$10^{-1}$,$10^0$,$10^1$},
      title={$\lambda$},title style={yshift=-2mm, xshift=0mm} },
    xlabel={$\norm{\ndft \zb f - \zb g}_{2,\bw}^2$},
    ylabel={$\norm{\grad \zb f}_{1,2}$},
    x label style={yshift=.5em}, y label style={yshift=-1em},
      ]
      \addplot+[scatter,thick,mark size=2,color=gray,point meta={lg10(\thisrowno{0})}, scatter/use mapped color=
      {draw opacity=0,fill=mapped color}
      ] table[x index=1,y index=2, meta index=0] {heart1-noise0.05-rec-tv-lcurve.txt};
    \end{loglogaxis}
  \end{tikzpicture}
  \caption{L-curve}
\end{subfigure}%
\begin{minipage}[b]{.55\textwidth}\centering
\begin{subfigure}[b]{.33\textwidth}\centering
  \begin{tikzpicture}
    \begin{axis}[
      width=\smallmodelwidth,
      height=\smallmodelwidth,
      enlargelimits=false,
      scale only axis,
      axis on top,
      ticks=none,
      ]
      \addplot graphics [
      xmin=\xmin, xmax=\xmax,  ymin=\xmin, ymax=\xmax,
      ] {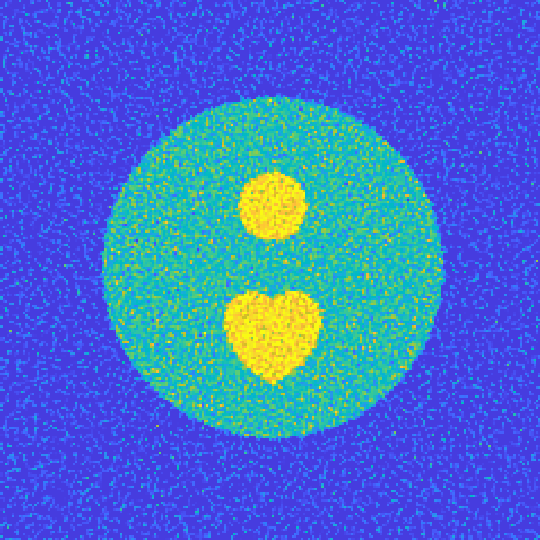};
    \end{axis}
  \end{tikzpicture}
  \caption{ $\lambda=0.001$}
\end{subfigure}%
\begin{subfigure}[b]{.33\textwidth}\centering
\begin{tikzpicture}
  \begin{axis}[
    width=\smallmodelwidth,
    height=\smallmodelwidth,
    enlargelimits=false,
    scale only axis,
    axis on top,
    ticks=none,
    ]
    \addplot graphics [
    xmin=\xmin, xmax=\xmax,  ymin=\xmin, ymax=\xmax,
    ] {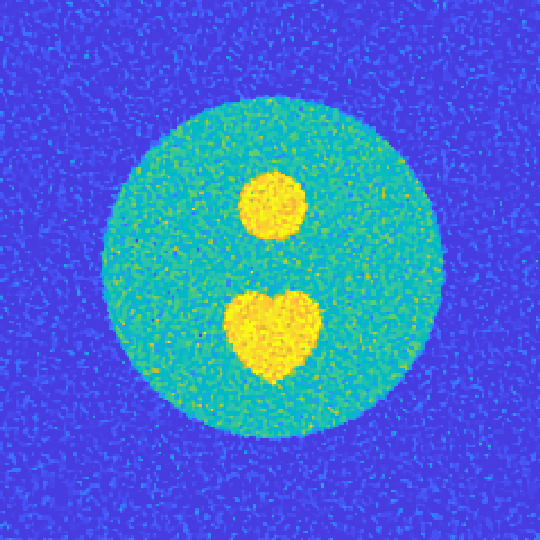};
  \end{axis}
\end{tikzpicture}
\caption{ $\lambda=0.01$}
\end{subfigure}%
\begin{subfigure}[b]{.33\textwidth}\centering
\begin{tikzpicture}
  \begin{axis}[
    width=\smallmodelwidth,
    height=\smallmodelwidth,
    enlargelimits=false,
    scale only axis,
    axis on top,
    ticks=none,
    ]
    \addplot graphics [
    xmin=\xmin, xmax=\xmax,  ymin=\xmin, ymax=\xmax,
    ] {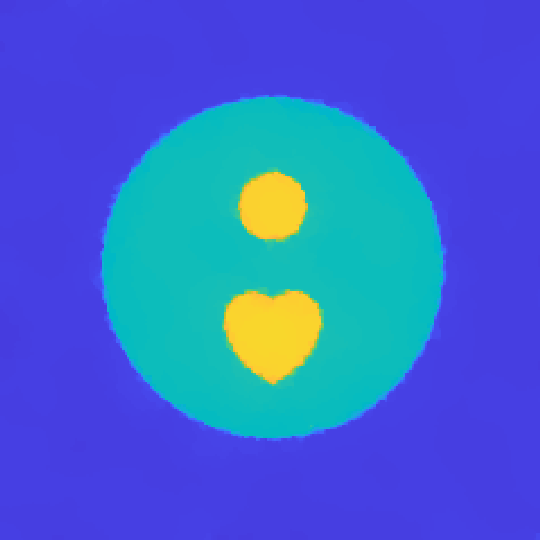};
  \end{axis}
\end{tikzpicture}
\caption{ $\lambda=0.1$}
\end{subfigure}\\
\begin{subfigure}[b]{.33\textwidth}\centering
\begin{tikzpicture}
  \begin{axis}[
    width=\smallmodelwidth,
    height=\smallmodelwidth,
    enlargelimits=false,
    scale only axis,
    axis on top,
    ticks=none,
    ]
    \addplot graphics [
    xmin=\xmin, xmax=\xmax,  ymin=\xmin, ymax=\xmax,
    ] {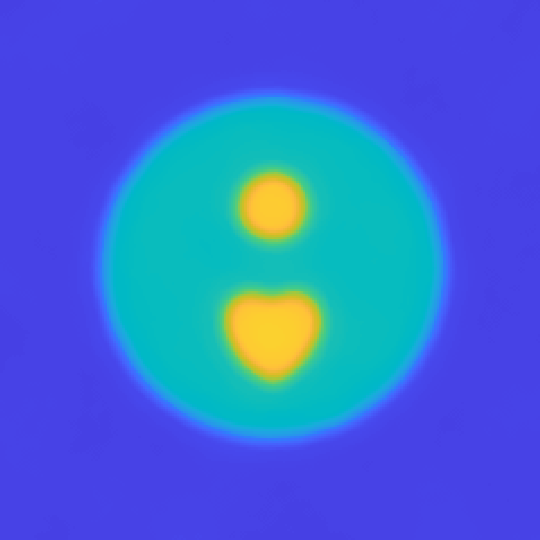};
  \end{axis}
\end{tikzpicture}
\caption{ $\lambda=1$}
\end{subfigure}
\begin{subfigure}[b]{.33\textwidth}\centering
\begin{tikzpicture}
  \begin{axis}[
    width=\smallmodelwidth,
    height=\smallmodelwidth,
    enlargelimits=false,
    scale only axis,
    axis on top,
    ticks=none,
    ]
    \addplot graphics [
    xmin=\xmin, xmax=\xmax,  ymin=\xmin, ymax=\xmax,
    ] {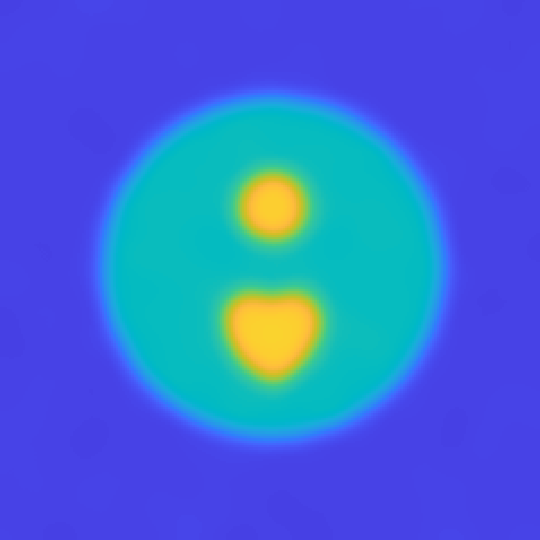};
  \end{axis}
\end{tikzpicture}
\caption{ $\lambda=10$}
\end{subfigure}
\end{minipage}
  \caption{L-curve for the PD with TV in \autoref{fig:heart1-noise-tv} for $\lambda\in[10^{-3},10]$.
    \label{fig:heart1-noise-lcurve}
  }
\end{figure}

Instead of choosing an appropriate $\lambda$ for PD with TV by hand as
in the above simulations, different choice strategies may be applied
\cite[Sect.\ 4]{EnHaNe96}.  We obtained promising results with the
L-curve method \cite{Han92}, which does not require a priori knowledge
of the noise level.  This method consists of choosing $\lambda$ with
respect to the ``tip'' of an L-shaped curve in a bilogarithmic plot of
the residual $\norm{\ndft \zb f_\lambda - \zb g}_{2,\bw}^2$ versus the
regularization $\norm{\grad \zb f_\lambda}_{1,2}$ for varying $\lambda$, where
$\zb f_\lambda$ denotes the reconstruction for a specific $\lambda$,
and where $\zb g$ is computed form the data $\bv^{\mathrm{sca}}_m$ as in
\autoref{alg:bw}.  An exemplary L-curve is displayed in
\autoref{fig:heart1-noise-lcurve}, where we can spot the ``tip''
around $\lambda \approx 0.1$.  A second, much weaker kink at the top
left around $\lambda\approx0.01$ can be easily ruled out as the
corresponding reconstruction $\zb f_\lambda$ looks very noisy.

\paragraph{Three-dimensional function.}

\begin{figure}[htp]\centering
  \begin{subfigure}[t]{.3\textwidth}\centering
    \begin{tikzpicture}
      \begin{axis}[
        width=.9*\modelwidth,
        height=.9*\modelwidth,
        enlargelimits=false,
        scale only axis,
        axis on top,
        colorbar,colorbar style={
          width=.15cm, xshift=-0.9em, 
          point meta min=0.6107,point meta max=3.4494,
        },
        ]
        \addplot graphics [
        xmin=-40, xmax=40,
        ymin=-40, ymax=40,
        ] {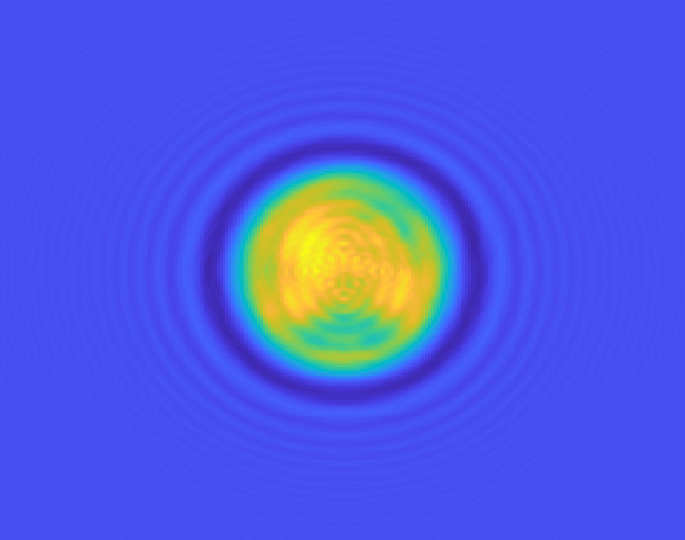};
      \end{axis}
    \end{tikzpicture}
    \caption{Exact data.}
  \end{subfigure}\hfill
  \begin{subfigure}[t]{.3\textwidth}\centering
    \begin{tikzpicture}
      \begin{axis}[
        width=.9*\modelwidth,
        height=.9*\modelwidth,
        enlargelimits=false,
        scale only axis,
        axis on top,
        colorbar,colorbar style={
          width=.15cm, xshift=-0.9em, 
          point meta min=0.5039,point meta max=3.5199,
        },
        ]
        \addplot graphics [
        xmin=-40, xmax=40,
        ymin=-40, ymax=40,
        ] {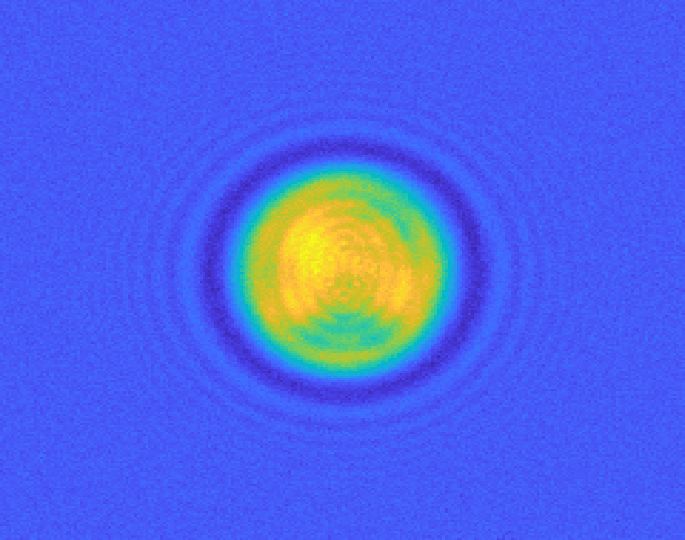};
      \end{axis}
    \end{tikzpicture}
    \caption{Data with 5\,\% Gaussian noise.}
  \end{subfigure}\hfill
  \begin{subfigure}[t]{.38\textwidth}\centering
  \includegraphics[width=\modelwidth]{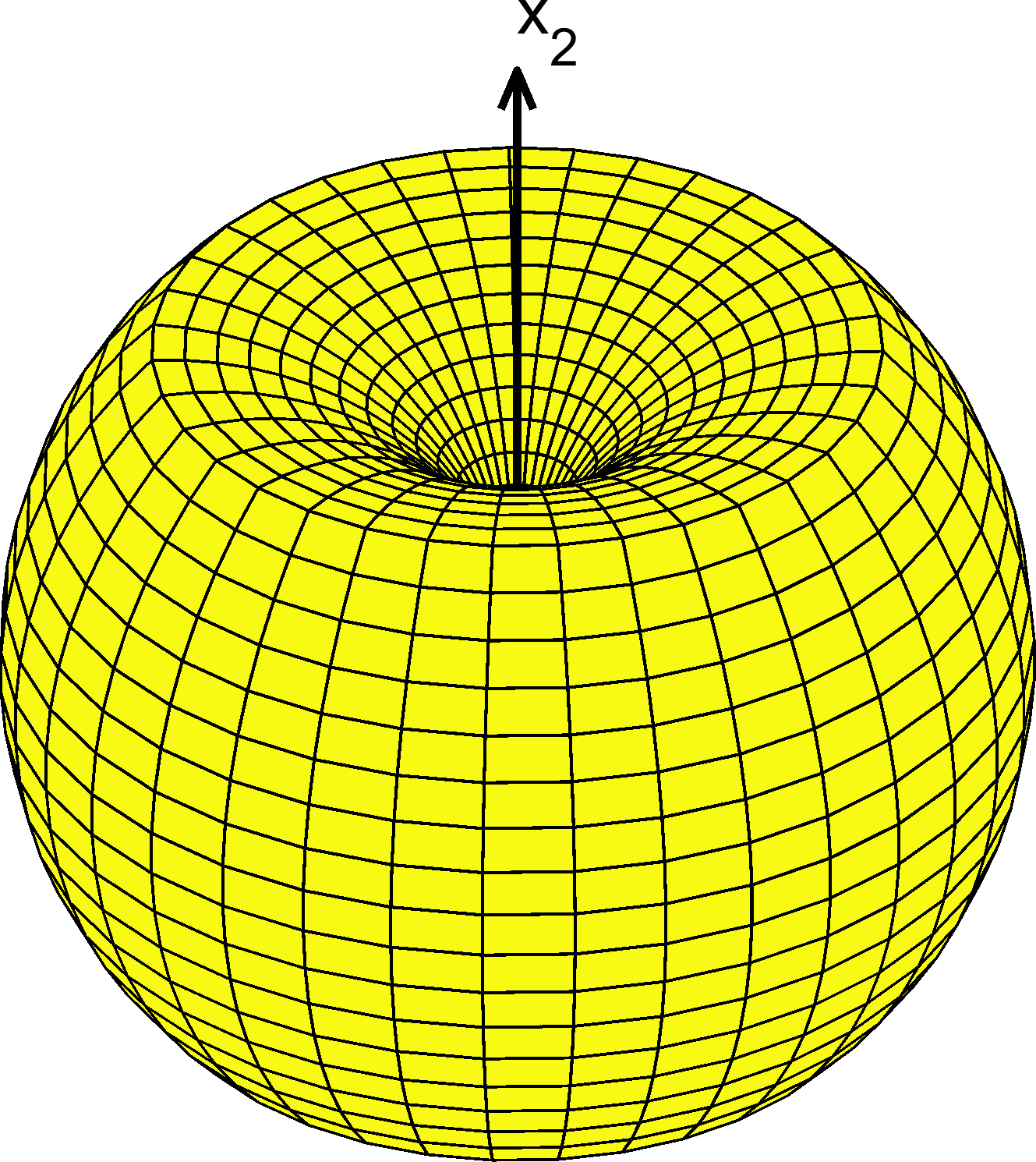}
  \caption{Fourier space coverage $\mathcal Y$ for a full turn around the $x_2$ axis. \label{fig:coverage3}}
  \end{subfigure}
  \caption{Data $\abs{\utot_t(\cdot,\rM)}$ for the function of \autoref{fig:cell3d-f} at one fixed time step $t$
    and the Fourier coverage $\mathcal Y$ in 3D. \label{fig:cell3d-sinogram}}
\end{figure}

\begin{figure}[ht]
  \pgfmathsetmacro{\cmin} {-0.05} \pgfmathsetmacro{\cmax} {1.05}
  \begin{subfigure}[t]{.33\textwidth}\centering
    \begin{tikzpicture}
      \begin{axis}[
        width=\modelwidth,
        height=\modelwidth,
        enlargelimits=false,
        scale only axis,
        axis on top,
        colorbar,colorbar style={
          width=.15cm, xshift=-0.9em, 
          point meta min=\cmin,point meta max=\cmax,
        },
        ]
        \addplot graphics [
        xmin=\xmin, xmax=\xmax, ymin=\xmin, ymax=\xmax,
        ] {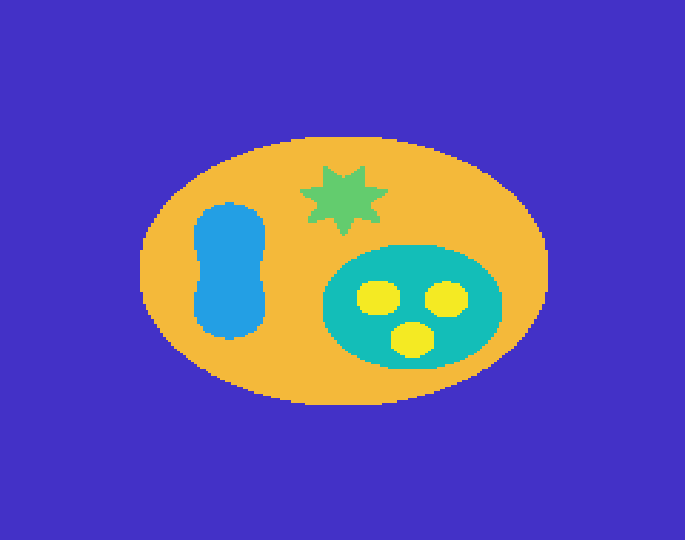};
      \end{axis}
    \end{tikzpicture}
    \caption{Ground truth $\bm f$ \label{fig:cell3d-f}}
  \end{subfigure}\hfill
  \begin{subfigure}[t]{.33\textwidth}\centering
    \begin{tikzpicture}
      \begin{axis}[
        width=\modelwidth,
        height=\modelwidth,
        enlargelimits=false,
        scale only axis,
        axis on top,
        colorbar,colorbar style={
          width=.15cm, xshift=-0.9em, 
          point meta min=\cmin,point meta max=\cmax,
        },
        ]
        \addplot graphics [
        xmin=\xmin, xmax=\xmax, ymin=\xmin, ymax=\xmax,
        ] {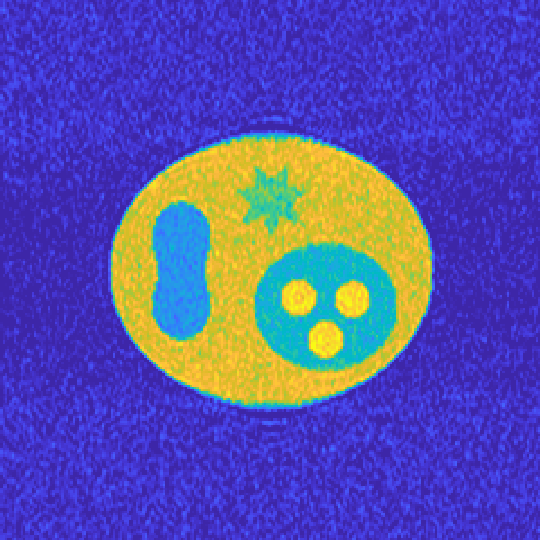};
      \end{axis}
    \end{tikzpicture}
    \caption{BP reconstruction\\PSNR 24.53, SSIM 0.178 \label{fig:cell3d-rec-bp}}
  \end{subfigure}\hfill
  \begin{subfigure}[t]{.33\textwidth}\centering
    \begin{tikzpicture}
      \begin{axis}[
        width=\modelwidth,
        height=\modelwidth,
        enlargelimits=false,
        scale only axis,
        axis on top,
        colorbar,colorbar style={
          width=.15cm, xshift=-0.9em, 
          point meta min=\cmin,point meta max=\cmax,
        },
        ]
        \addplot graphics [
        xmin=\xmin, xmax=\xmax, ymin=\xmin, ymax=\xmax,
        ] {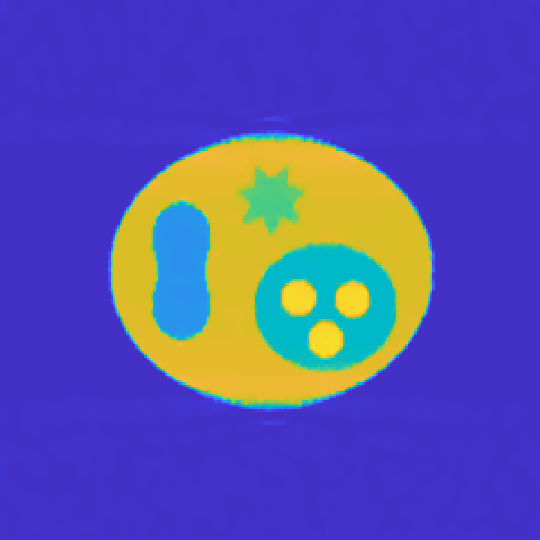};
      \end{axis}
    \end{tikzpicture}
    \caption{BP and TVd ($\lambda=0.05$)\\ PSNR 32.51, SSIM 0.968}
  \end{subfigure}

  \begin{subfigure}[t]{.33\textwidth}\centering
    \begin{tikzpicture}
      \begin{axis}[
        width=\modelwidth,
        height=\modelwidth,
        enlargelimits=false,
        scale only axis,
        axis on top,
        colorbar,colorbar style={
          width=.15cm, xshift=-0.9em, 
          point meta min=\cmin,point meta max=\cmax,
        },
        ]
        \addplot graphics [
        xmin=\xmin, xmax=\xmax, ymin=\xmin, ymax=\xmax,
        ] {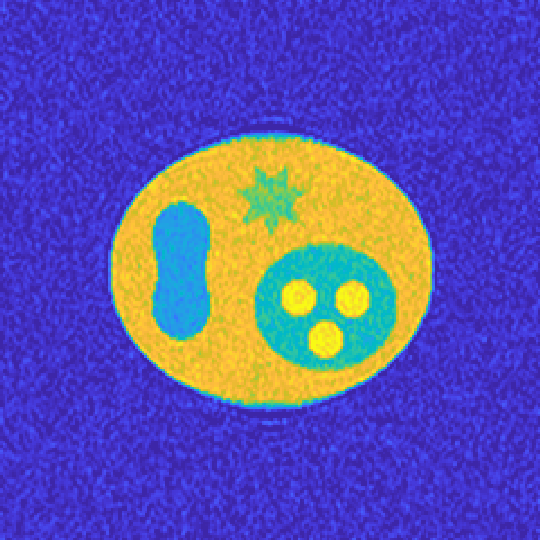};
      \end{axis}
    \end{tikzpicture}
    \caption{CG Reconstruction\\PSNR 26.74, SSIM 0.309 \label{fig:cell3d-rec-cg}}
  \end{subfigure}\hfill
  \begin{subfigure}[t]{.33\textwidth}\centering
    \begin{tikzpicture}
      \begin{axis}[
        width=\modelwidth,
        height=\modelwidth,
        enlargelimits=false,
        scale only axis,
        axis on top,
        colorbar,colorbar style={
          width=.15cm, xshift=-0.9em, 
          point meta min=\cmin,point meta max=\cmax,
        },
        ]
        \addplot graphics [
        xmin=\xmin, xmax=\xmax, ymin=\xmin, ymax=\xmax,
        ] {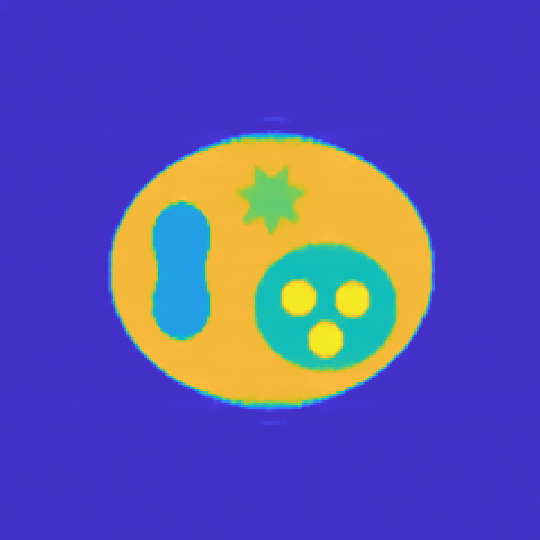};
      \end{axis}
    \end{tikzpicture}
    \caption{CG and TVd ($\lambda=0.02$)\\ PSNR 33.53, SSIM 0.965}
  \end{subfigure}\hfill
  \begin{subfigure}[t]{.33\textwidth}\centering
    \begin{tikzpicture}
      \begin{axis}[
        width=\modelwidth,
        height=\modelwidth,
        enlargelimits=false,
        scale only axis,
        axis on top,
        colorbar,colorbar style={
          width=.15cm, xshift=-0.9em, 
          point meta min=\cmin,point meta max=\cmax,
        },
        ]
        \addplot graphics [
        xmin=\xmin, xmax=\xmax, ymin=\xmin, ymax=\xmax,
        ] {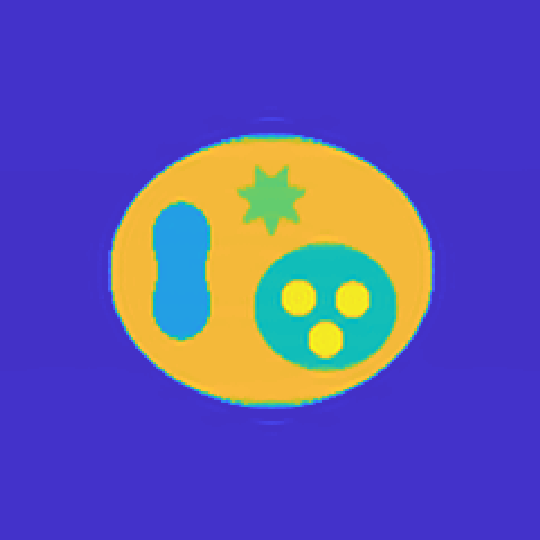};
      \end{axis}
    \end{tikzpicture}
    \caption{PD with TV ($\lambda=0.02$)\\PSNR 33.58, SSIM 0.736}
  \end{subfigure}%
  \caption{Slice $x_3=0.35$ of the 3D reconstruction from data with 5\,\% Gaussian noise using BP in \eqref{eq:bp-d}, CG in Alg.~\ref{alg:cg}, and PD
  with TV in Alg.~\ref{alg:pd-back}.  
  \label{fig:cell3d-noise}}
\end{figure}

Our three-dimensional phantom contains a variety of shapes with smooth
and non-smooth boundaries, see \autoref{fig:cell3d-f}.  The
discretization parameters are the same as before.  The object makes a
full turn around the $x_2$ axis with the rotation angle $t\in[0,2\pi]$.  The
measurements with respect to one illumination angle is shown in
\autoref{fig:cell3d-sinogram}.
The respective Fourier space
coverage~$\mathcal Y$ in \eqref{eq:Y} has the form of a horn torus
and is depicted in \autoref{fig:coverage3}, cf.\
\cite{KirQueRitSchSet21}.  The so-called \emph{missing cone} around the
$x_2$ axis in the Fourier space coverage makes the reconstruction more
difficult and promotes certain artifacts, like horizontal stripes at
the top and bottom in \autoref{fig:cell3d-noise}.

\begin{table}[htp]\centering
  \begin{tabular}{c c c c c c}
    & BP & BP\,\&\,TVd & CG & CG\,\&\,TVd & PD
    \\\hline
    PSNR
    & 29.07 & 32.79 & 33.36 & 33.82 & 34.36
    \\\hline SSIM
    & 0.614 & 0.987 & 0.955 & 0.988 & 0.957
    \\\hline Time (sec)
    & 4 & 31 & 79 & 31 & 1395
  \end{tabular}
  \caption{Error for the reconstruction of the function of \autoref{fig:cell3d-f} with exact data.
    \label{tab:cell3d}}
\end{table}

The errors for exact data and computation times are presented in \autoref{tab:cell3d}.
The reconstructions with noisy data are shown in \autoref{fig:cell3d-noise}.  
The CG
reconstruction is again of better quality than BP.  Both, however,
contain visible noise, 
which is compensated in the
additional TV denoising performed afterwards, where we make
$J_{\mathrm{PD}}=20$ iterations.  Performing more iterations does not
benefit the image quality significantly.  The PD method with TV for
the full inverse problem produces also a very good 3D reconstruction.
However, the latter requires significantly more time since the
primal-dual iteration converges rather slowly and requires
$J_{\mathrm{PD}} = 100$ iterations to obtain a satisfying result.
Furthermore, we notice that TV regularization may lead to staircasing \cite{Jal16}
if the image is not piecewise constant.
This phenomenon might be prevented by using higher-order generalizations of the total variation, see \cite{BreHol20}.

\section{Phase retrieval in ODT}
\label{sec:phase-retrieval-odt}

Different from the numerical experiments in
\autoref{sec:numer-backw-transf}, the phase of the complex diffraction
patterns $\utot_t$ is usually unknown in practice.  In other words, only
the intensity $\abs{\utot_t(\cdot,\rM)}$ is recorded.  Due to the
phase loss, the already ill-posed recovery problem is further
hampered.

The numerical retrieval of an unknown phase traces back to the famous
Gerchberg--Saxton algorithm \cite{GerSax72}, which was originally
invented to capture a complex-valued object from its intensities in
spatial and Fourier domain and was modified by Fienup
\cite{Fie78} to recover an object from the modulus of its Fourier
transform.  The basic idea of Fienup's error-reduction algorithm is to
enforce known object constraints in object domain and known magnitude
constraints given by the data in Fourier domain alternatingly.  The
error-reduction algorithm is a specific version of more general
input-output algorithms, which may be read as \autoref{alg:err-red},
where $\bm D$ is a map modeling the measurement process, and
where the sign function
\begin{equation*}
  \sgn(z) \coloneqq
  \begin{cases}
    z / |z|, & z \in\C\setminus \{0\}, \\
    1, & z = 0,
  \end{cases}
\end{equation*}
is applied element by element.  The error-reduction algorithm and
modifications of it are nowadays still well-established approaches to
solve phase retrieval problems in practice.

\begin{algorithm}[hbt]
  \KwIn{$\bm d = \abs{\bm D(\bm f)}$.}
  
  Initialize $\bm g^{(0)} \coloneqq \bm d$\;
  
  \For{$j=0,1,2,\dots,J_{\mathrm{IO}}$}{
    $\bm f^{(j)} \coloneqq \bm D^{-1} \bm g^{(j)}$, see \autoref{alg:bw}\;
    Set $\bm f^{(j+1/2)}$ by applying object domain constraints on $\bm f^{(j)}$ via
    \eqref{eq:f-projection} or \eqref{eq:hio-update}\;
    $\bm g^{(j+1/2)} \coloneqq \bm D \bm f^{(j+1/2)}$, see \autoref{alg:fw}\;
    $\bm g^{(j+1)} \coloneqq \bm d \sgn (\bm g^{(j+1/2)})$\;
  }
  \KwOut{Approximate scattering potential $\bm f^{(J_{\mathrm{IO}})}$.}
  \caption{General input-output algorithm} 
  \label{alg:err-red}
\end{algorithm}

For our specific application, the given data $\bm d = |\bm D(\bm f)|$
corresponds to the measurement process modeled via
$\bm D^{\mathrm{tot}}$ in \eqref{eq:Dtot}.  Since
$\bm D^{\mathrm{tot}}$ is not invertible in general, we approximately
compute $\bm f^{(j)}$ via \autoref{alg:bw}.  The object domain
constraints of $\bm f$ consist in the non-negativity and support
constraint of the investigated physical object.  Assuming that the
support of $f$ is constrained in the ball of radius $\rs$, we may enforce
them by setting
\begin{equation} \label{eq:f-projection}
  \tilde f_\bk^{(j)}  \coloneqq
  \begin{cases}
    \max\{f^{(j)}_{\bk} , 0 \}, & \| \bx_\bk \|_2 \le \rs, \\
    0, & \| \bx_\bk \|_2 > \rs,
  \end{cases}
\end{equation}
for $\bk\in\mathcal I_K^d$.
In the following, we refer to \autoref{alg:err-red} with the input rule
$\bm f^{(j+1/2)} \coloneqq \tilde{\bm f}^{(j)}$ as the \emph{error-reduction algorithm (for
  ODT)}.

Instead of enforcing the object domain constrains strictly, the input
of the next forward transform may be modified such that the output of
the subsequent backward transform is pushed in direction of the
feasible object signals.  This procedure usually leads to a better
numerical convergence of \autoref{alg:err-red}.  A typical choice of
the next input $f^{(j+1/2)}$ is known as the \emph{hybrid input-output} (HIO)
by Fienup \cite{Fie82}.  For some fixed $\beta \in (0,1]$, the HIO
step in \autoref{alg:err-red} consists of
\begin{equation}
  \label{eq:hio-update}
  f^{(j+1/2)}_{\bk} \coloneqq
  \begin{cases}
    f^{(j)}_{\bk}, & \bk \in \Gamma^{(j)},\\
    f^{(j-1/2)}_{\bk} - \beta (f^{(j)}_{\bk} - \tilde f^{(j)}_{\bk}),
    & \bk \in \mathcal I_K^d \setminus\Gamma^{(j)},
  \end{cases}
\end{equation}
where
$\Gamma^{(j)} \coloneqq \big\{\bk\in\mathcal I_K^d : f^{(j)}_{\bk} = \tilde
f^{(j)}_{\bk} \big\}$ consists of all indices $\bk$ where
$f^{(j)}_\bk$ satisfies the object domain constraints.  The next input
$\bm f^{(j+1/2)}$ is thus composed of the input $\bm f^{(j-1/2)}$ and
output $\bm f^{(j)}$ of the last forward-backward transform. We
henceforth refer to \autoref{alg:err-red} with input rule
\eqref{eq:hio-update} as the \emph{HIO algorithm (for ODT)}.

The numerical performance of the error-reduction and the hybrid
input-output algorithm may be increased by denoising the output of the
forward-backward transform using the total variation.  This idea has
been successfully applied for example in Fourier and Fresnel diffraction
imaging \cite{GauMohKha15,GauKha19}.  If the inverse of the NDFT
during the backward transform is computed by \autoref{alg:pd-back}, we
automatically obtain a TV regularized output; so the backward
transform corresponding to \eqref{eq:var-prob} goes hand in hand with
the TV heuristic for the input-output phase retrieval methods.

\section{Numerics of phase retrieval}
\label{sec:numer-pr}

In this section, we compare the performance of the CG method and the
PD method in the context of phase retrieval in ODT, where we rely on
the hybrid input-output algorithm in \autoref{sec:phase-retrieval-odt}.
In the second subsection, we compare our method with an existing one from literature, 
which only works well for a large distance $\rM$ of the measurement plane.

\subsection{Full diffraction model}
\label{sec:full-diffr-model}

We apply \autoref{alg:err-red} for the diffraction
model $\zb D^\mathrm{tot}$.  For the inversion of
$\zb D^\mathrm{tot}$ required in each step of the input-output
algorithm, we utilize either the CG method or the PD method with TV.
We do not show the results with the BP method here, because it  always performed worse than the CG.
Since CG and PD are inside the loop of the input-output algorithm, the
number of inner iterations can be reduced.  In particular, the PD
method in \autoref{alg:pd-back} may be restarted in each step of the
outer loop with the previous dual variable.  Instead of restarting,
the backtracking method may also be resumed, which significantly
reduces the run time of the overall algorithm.  Since the CG-based phase
retrieval turns out to be much faster than the PD-based method, we
utilize the CG phase reconstruction as starting point for phase
retrieval with PD.

\begin{figure}[t]\centering
  \pgfmathsetmacro{\cmin} {-0.05} \pgfmathsetmacro{\cmax} {0.55}
  \begin{subfigure}[t]{.33\textwidth}
    \begin{tikzpicture}
      \begin{axis}[
        width=\modelwidth,
        height=\modelwidth,
        enlargelimits=false,
        scale only axis,
        axis on top,
        colorbar,colorbar style={
          width=.15cm, xshift=-0.9em, 
          point meta min=\cmin,point meta max=\cmax,
        },
        ]
        \addplot graphics [
        xmin=\xmin, xmax=\xmax,  ymin=\xmin, ymax=\xmax,
        ] {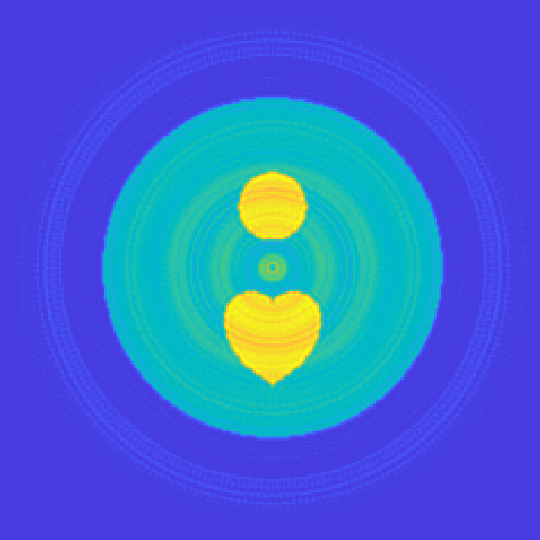};
      \end{axis}
    \end{tikzpicture}
    \caption{HIO/CG reconstruction\\  $J_{\mathrm{IO}}=10$, $J_{\mathrm{CG}}=5$\\PSNR 34.44, SSIM 0.821 \\0.6\,sec}
  \end{subfigure}\hspace{-6pt}
  \begin{subfigure}[t]{.33\textwidth}\centering
    \begin{tikzpicture}
      \begin{axis}[
        width=\modelwidth,
        height=\modelwidth,
        enlargelimits=false,
        scale only axis,
        axis on top,
        colorbar,colorbar style={
          width=.15cm, xshift=-0.9em, 
          point meta min=\cmin,point meta max=\cmax,
        },
        ]
        \addplot graphics [
        xmin=\xmin, xmax=\xmax,  ymin=\xmin, ymax=\xmax,
        ] {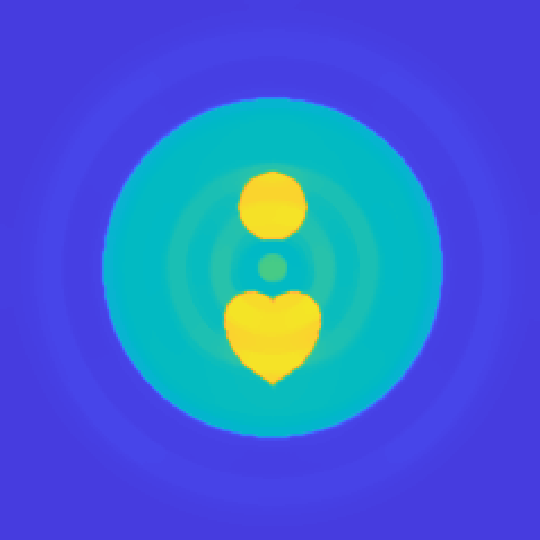};
      \end{axis}
    \end{tikzpicture}
  \caption{HIO/CG and TVd \\ $\lambda=0.05$, $J_{\mathrm{PD}}=50$\\PSNR 35.42, SSIM 0.831 \\0.1\,sec}
  \end{subfigure}\hspace{-6pt}
  \begin{subfigure}[t]{.33\textwidth}\centering
    \begin{tikzpicture}
      \begin{axis}[
        width=\modelwidth,
        height=\modelwidth,
        enlargelimits=false,
        scale only axis,
        axis on top,
        colorbar,colorbar style={
          width=.15cm, xshift=-0.9em, 
          point meta min=\cmin,point meta max=\cmax,
        },
        ]
        \addplot graphics [
        xmin=\xmin, xmax=\xmax,  ymin=\xmin, ymax=\xmax,
        ] {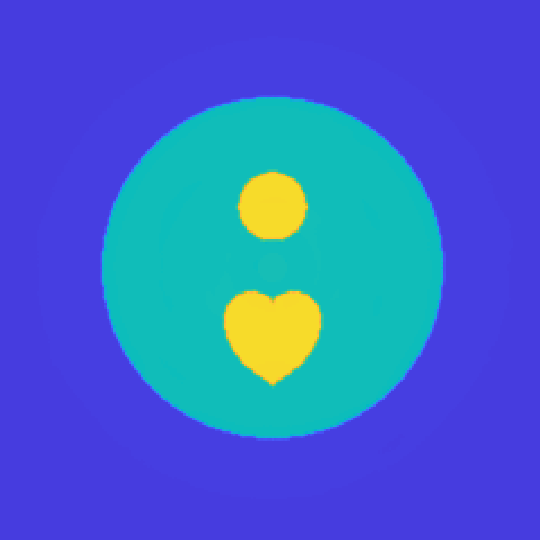};
      \end{axis}
    \end{tikzpicture}
    \caption{HIO/PD with TV\\$\lambda{=}0.01$, $J_{\mathrm{IO}}{=}20$, $J_{\mathrm{PD}}{=}5$\\PSNR 41.32, SSIM 0.981 \\ 13\,sec}
  \end{subfigure}
  \caption{HIO phase retrieval of $|\bm D^{\mathrm{tot}} \bm f|$ using CG in Alg.~\ref{alg:cg} and PD
  with TV in Alg.~\ref{alg:pd-back}.  The HIO/CG result is 
  TV-denoised.
  \label{fig:heart1-pr}}
\end{figure}

\begin{figure}[t]\centering
  \pgfmathsetmacro{\cmin} {-0.05} \pgfmathsetmacro{\cmax} {0.55}
  \begin{subfigure}[t]{.33\textwidth}
    \begin{tikzpicture}
      \begin{axis}[
        width=\modelwidth,
        height=\modelwidth,
        enlargelimits=false,
        scale only axis,
        axis on top,
        colorbar,colorbar style={
          width=.15cm, xshift=-0.9em, 
          point meta min=\cmin,point meta max=\cmax,
        },
        ]
        \addplot graphics [
        xmin=\xmin, xmax=\xmax,  ymin=\xmin, ymax=\xmax,
        ] {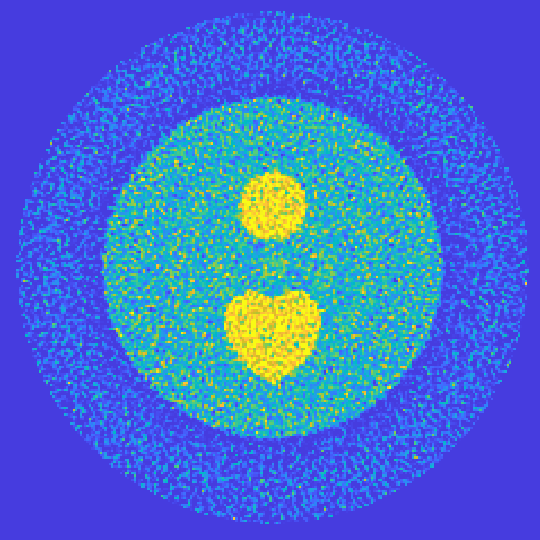};
      \end{axis}
    \end{tikzpicture}
    \caption{HIO/CG reconstruction\\  $J_{\mathrm{IO}}=10$, $J_{\mathrm{CG}}=5$\\PSNR 22.49, SSIM 0.354}
  \end{subfigure}\hspace{-8pt}
  \begin{subfigure}[t]{.33\textwidth}\centering
    \begin{tikzpicture}
      \begin{axis}[
        width=\modelwidth,
        height=\modelwidth,
        enlargelimits=false,
        scale only axis,
        axis on top,
        colorbar,colorbar style={
          width=.15cm, xshift=-0.9em, 
          point meta min=\cmin,point meta max=\cmax,
        },
        ]
        \addplot graphics [
        xmin=\xmin, xmax=\xmax,  ymin=\xmin, ymax=\xmax,
        ] {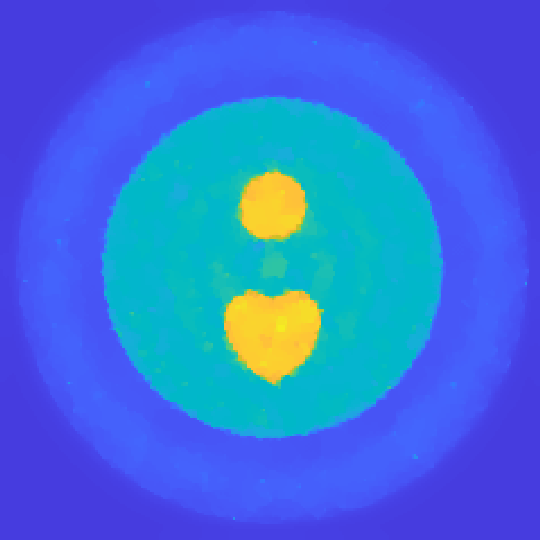};
      \end{axis}
    \end{tikzpicture}
    \caption{HIO/CG and TVd\\ $\lambda=0.1$, $J_{\mathrm{PD}} = 50$\\PSNR 28.38, SSIM 0.572}
  \end{subfigure}\hspace{-8pt}
  \begin{subfigure}[t]{.33\textwidth}\centering
    \begin{tikzpicture}
      \begin{axis}[
        width=\modelwidth,
        height=\modelwidth,
        enlargelimits=false,
        scale only axis,
        axis on top,
        colorbar,colorbar style={
          width=.15cm, xshift=-0.9em, 
          point meta min=\cmin,point meta max=\cmax,
        },
        ]
        \addplot graphics [
        xmin=\xmin, xmax=\xmax,  ymin=\xmin, ymax=\xmax,
        ] {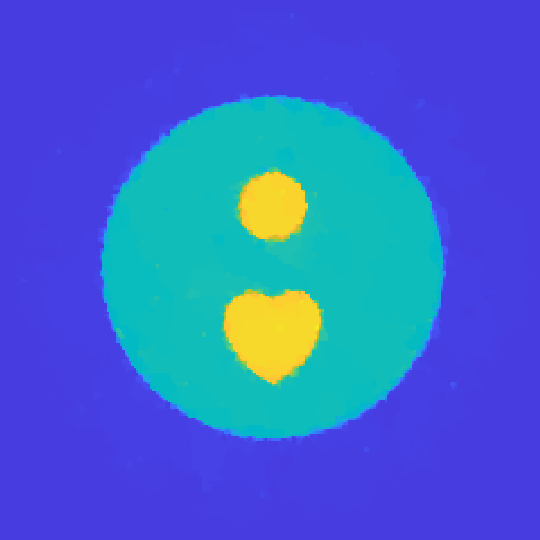};
      \end{axis}
    \end{tikzpicture}
    \caption{HIO/PD with TV\\
      $\lambda={0}.05$, $J_{\mathrm{IO}}{=}50$,
        $J_{\mathrm{PD}}{=}10$\\PSNR 37.27, SSIM
      0.936\label{fig:heart1-noise-pr:PD}}
  \end{subfigure}
  \caption{HIO phase retrieval of $|\bm D^{\mathrm{tot}} \bm f|$ with
    5\% Gaussian noise using CG in Alg.~\ref{alg:cg} and PD
  with TV in Alg.~\ref{alg:pd-back}.  The HIO/CG result is 
  TV-denoised.
  \label{fig:heart1-noise-pr}}
\end{figure}

We employ the problem setting from \autoref{sec:numer-backw-transf}.
It is especially relevant for phase retrieval that the images are
padded with zeros to represent the support constraints
\eqref{eq:f-projection}.  We impose the support constraint $\rs = 40$
in \eqref{eq:f-projection}, which is larger than the actual support
radius of the ground truth function of about $26$.  The results are
shown in \autoref{fig:heart1-pr},
where the HIO is applied with the parameter $\beta=0.7$.
Especially for PD with TV, the
result nearly coincides with the ground truth; so the missing phase
does not decrease the image quality too much.  With the CG method,
there are some visible artifacts.  However, executing more CG
iterations may decrease the image quality because of the ill-posedness of
the problem.
The HIO/CG method coupled with TV denoising is outperformed by the HIO/PD method since the latter employs TV in each step of the outer iteration,
while the former does it only once.

Adding 5\,\% additive Gaussian noise to the intensity data $\abs{\zb D^{\mathrm{tot}}\zb f}$ does
not decrease the reconstruction quality of HIO with TV-regularized PD
inversion much, whereas the HIO reconstruction with CG becomes worse,
see \autoref{fig:heart1-noise-pr}.  Raising the number of HIO
iterations with CG results in even noisier reconstructions.
For the PD method, the small number of inner iterations $J_{\mathrm{PD}}=10$ is sufficient because we restart the inner loop with the dual variables, parameters and initial solution from the previous outer step.
Without this resuming, we need considerably more iterations, 
and even then the results are slightly worse,
see \autoref{tab:heart1-noise-pr}.

\begin{table}[ht]\centering
  \begin{tabular}{r | c c c c | c}
    & \multicolumn{4}{c|}{cold start with $J_{\mathrm{PD}}=$} & warm start
    \\
    $J_{\mathrm{IO}}=50$ & 10 & 20 & 50 & 100
                                        & $J_{\mathrm{PD}}=10$
    \\\hline PSNR
    & 29.64 & 35.10 & 36.34 & 36.44 
                                        & 37.12
    \\ SSIM
    & 0.687 & 0.813 & 0.841 & 0.848 
                                        & 0.915
  \end{tabular}
  \caption{Quality of the HIO/PD phase retrieval for cold start of the PD depending on the number of inner iterations $J_{\mathrm{PD}}$
  and for the warm start (resuming with parameters, dual variable and initial solution of the previous step) 
  as in \autoref{fig:heart1-noise-pr:PD}.
  \label{tab:heart1-noise-pr}}
\end{table}

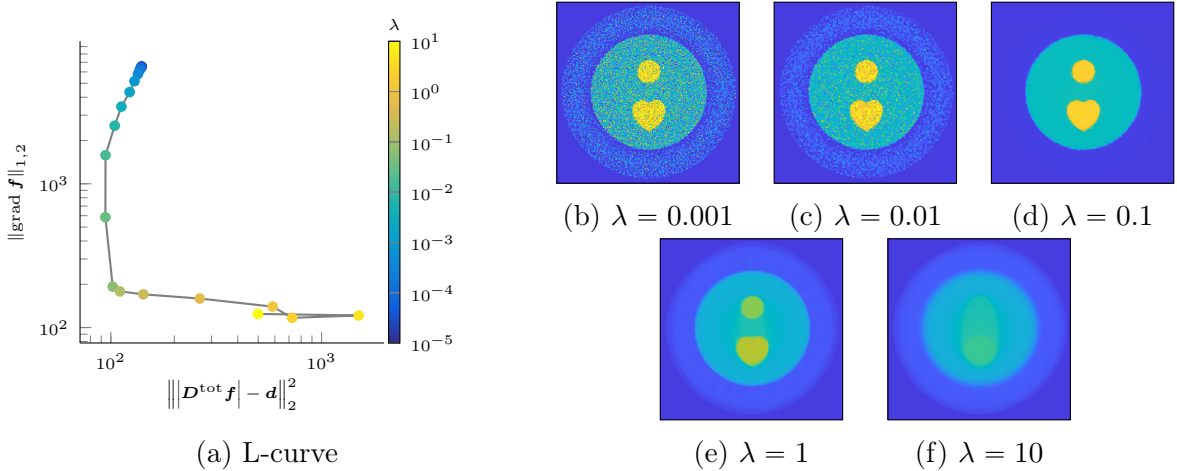
\begin{figure}[t]\centering
  \pgfmathsetmacro{\cmin} {-0.05} \pgfmathsetmacro{\cmax} {0.55}
  \begin{subfigure}[b]{.45\textwidth}
  \begin{tikzpicture}
    \pgfmathdeclarefunction{lg10}{1}{\pgfmathparse{ln(#1)/ln(10)}}
    \begin{loglogaxis}[width=4cm, height=4cm, enlarge y limits=0.1, axis x line*=bottom, axis y line*=left, legend pos=outer north east, legend style={cells={anchor=west}}, scale only axis,
      colorbar, colorbar style={width=.15cm, xshift=-0.6em,
        ytick={-5,-4,-3,-2,-1,0,1},yticklabels={$10^{-5}$,$10^{-4}$,$10^{-3}$,$10^{-2}$,$10^{-1}$,$10^0$,$10^1$},
        title={$\lambda$},title style={yshift=-2mm, xshift=0mm} },
      xlabel={$\norm{\abs{\zb D^{\mathrm{tot}}\zb f} - \zb d}_{2}^2$},
      ylabel={$\norm{\grad \zb f}_{1,2}$},
      x label style={yshift=.5em}, y label style={yshift=-1em},
      ]
      \addplot+[scatter,thick,color=gray,mark size=2,point meta={lg10(\thisrowno{0})}, scatter/use mapped color=
      {draw opacity=0,fill=mapped color}
      ] table[x index=1,y index=2, meta index=0] {heart1-noise0.05-pr-tv-lcurve.txt};
    \end{loglogaxis}
  \end{tikzpicture}
  \caption{L-curve}
  \end{subfigure}%
  \begin{minipage}[b]{.55\textwidth}\centering
  \begin{subfigure}[b]{.33\textwidth}\centering
    \begin{tikzpicture}
      \begin{axis}[
        width=\smallmodelwidth,
        height=\smallmodelwidth,
        enlargelimits=false,
        scale only axis,
        axis on top,
        ticks=none,
        ]
        \addplot graphics [
        xmin=\xmin, xmax=\xmax,  ymin=\xmin, ymax=\xmax,
        ] {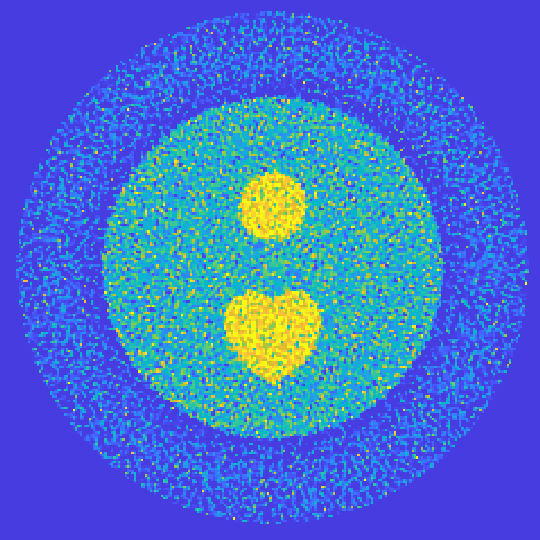};
      \end{axis}
    \end{tikzpicture}
    \caption{ $\lambda=0.001$}
  \end{subfigure}%
  \begin{subfigure}[b]{.33\textwidth}\centering
    \begin{tikzpicture}
      \begin{axis}[
        width=\smallmodelwidth,
        height=\smallmodelwidth,
        enlargelimits=false,
        scale only axis,
        axis on top,
        ticks=none,
        ]
        \addplot graphics [
        xmin=\xmin, xmax=\xmax,  ymin=\xmin, ymax=\xmax,
        ] {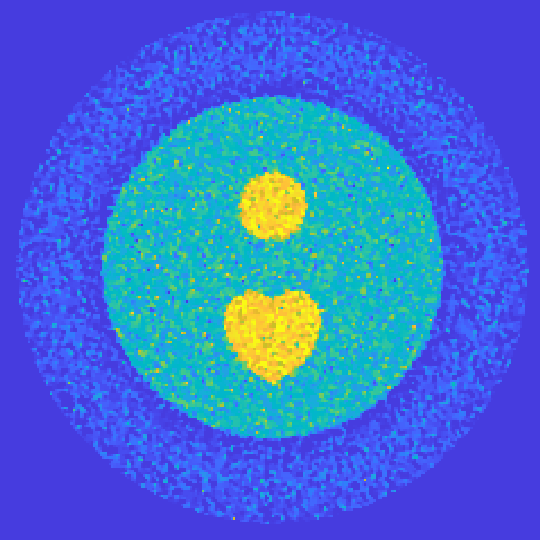};
      \end{axis}
    \end{tikzpicture}
    \caption{ $\lambda=0.01$}
  \end{subfigure}%
  \begin{subfigure}[b]{.33\textwidth}\centering
    \begin{tikzpicture}
      \begin{axis}[
        width=\smallmodelwidth,
        height=\smallmodelwidth,
        enlargelimits=false,
        scale only axis,
        axis on top,
        ticks=none,
        ]
        \addplot graphics [
        xmin=\xmin, xmax=\xmax,  ymin=\xmin, ymax=\xmax,
        ] {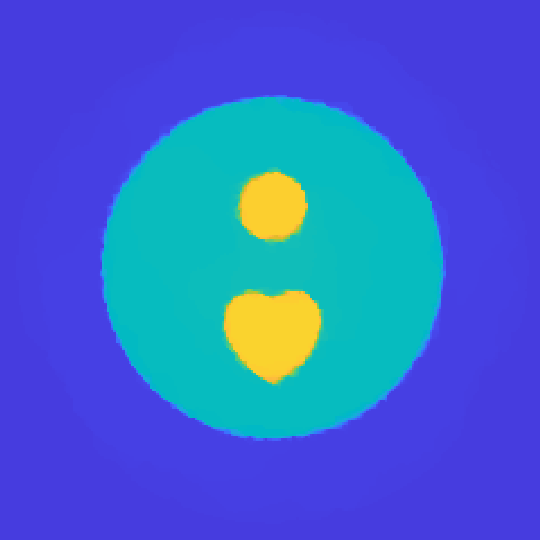};
      \end{axis}
    \end{tikzpicture}
    \caption{ $\lambda=0.1$}
  \end{subfigure}\\
  \begin{subfigure}[b]{.33\textwidth}\centering
    \begin{tikzpicture}
      \begin{axis}[
        width=\smallmodelwidth,
        height=\smallmodelwidth,
        enlargelimits=false,
        scale only axis,
        axis on top,
        ticks=none,
        ]
        \addplot graphics [
        xmin=\xmin, xmax=\xmax,  ymin=\xmin, ymax=\xmax,
        ] {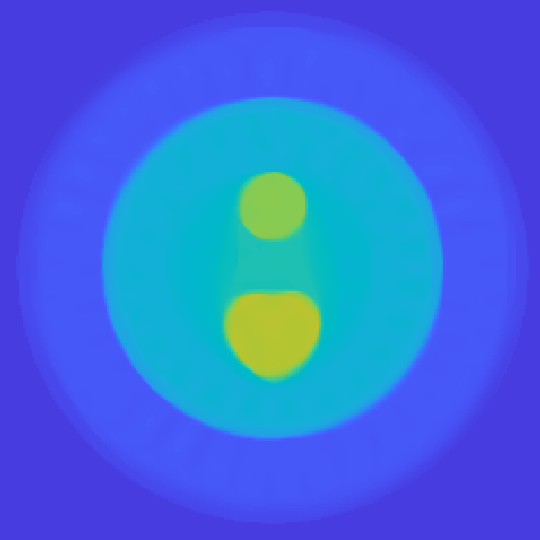};
      \end{axis}
    \end{tikzpicture}
    \caption{ $\lambda=1$}
  \end{subfigure}
  \begin{subfigure}[b]{.33\textwidth}\centering
    \begin{tikzpicture}
      \begin{axis}[
        width=\smallmodelwidth,
        height=\smallmodelwidth,
        enlargelimits=false,
        scale only axis,
        axis on top,
        ticks=none,
        ]
        \addplot graphics [
        xmin=\xmin, xmax=\xmax,  ymin=\xmin, ymax=\xmax,
        ] {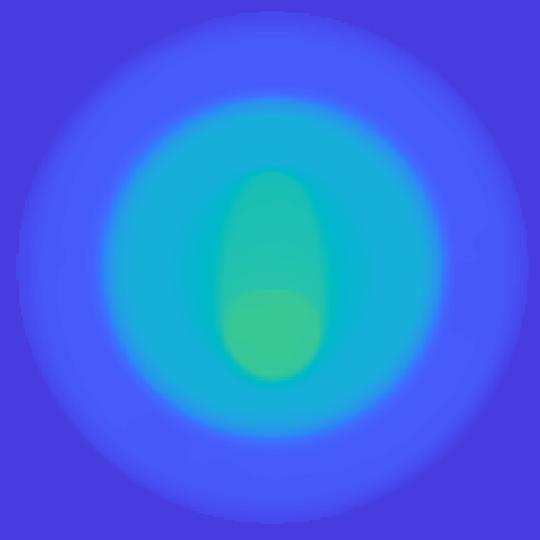};
      \end{axis}
    \end{tikzpicture}
    \caption{ $\lambda=10$}
  \end{subfigure}
\end{minipage}
\caption{L-curve for HIO/PD with TV in
  \autoref{fig:heart1-noise-pr:PD} for $\lambda\in[10^{-5},10]$.
    \label{fig:heart1-noise-pr-lcurve}
  }
\end{figure}

The L-curve method does not formally apply to this setting since the
regularization parameter $\lambda$ appears only in the PD method and
thus in the inner loop of \autoref{alg:err-red}.  However, we could
still achieve some reasonable results, see
\autoref{fig:heart1-noise-pr-lcurve}, where we plotted the residual
$\norm{\abs{\zb D^{\mathrm{tot}}\zb f} - \zb d}_{2}^2$ versus the
discrete TV semi-norm $\norm{\zb f}_{2,1}$.  Note that there is a
different behavior of the residual for very small and very large
$\lambda$, where the HIO algorithm fails to find the right phase.
However, these values may be excluded since no meaningful
reconstruction occurs.

\begin{figure}[t]\centering
  \pgfmathsetmacro{\cmin} {-0.05} \pgfmathsetmacro{\cmax} {1.05}
  \begin{subfigure}[t]{.33\textwidth}
    \begin{tikzpicture}
      \begin{axis}[
        width=\modelwidth,
        height=\modelwidth,
        enlargelimits=false,
        scale only axis,
        axis on top,
        colorbar,colorbar style={
          width=.15cm, xshift=-0.9em, 
          point meta min=\cmin,point meta max=\cmax,
        },
        ]
        \addplot graphics [
        xmin=\xmin, xmax=\xmax,  ymin=\xmin, ymax=\xmax,
        ] {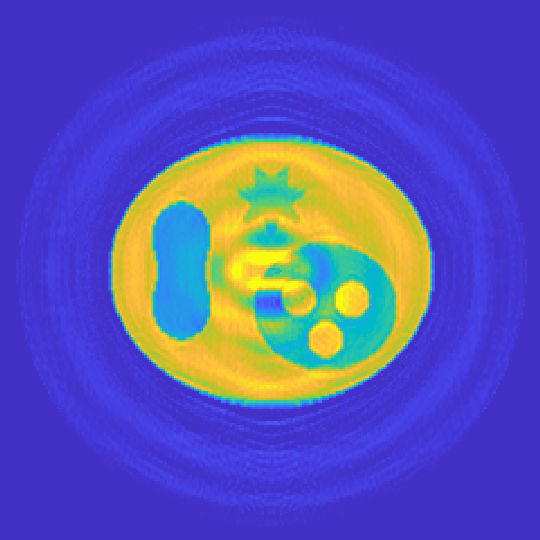};
      \end{axis}
    \end{tikzpicture}
    \caption{HIO/CG reconstruction\\  $J_{\mathrm{IO}}=10$, $J_{\mathrm{CG}}=5$\\PSNR 29.52, SSIM 0.713\\ 5\,min}
  \end{subfigure}\hspace{-6pt}
  \begin{subfigure}[t]{.33\textwidth}\centering
    \begin{tikzpicture}
      \begin{axis}[
        width=\modelwidth,
        height=\modelwidth,
        enlargelimits=false,
        scale only axis,
        axis on top,
        colorbar,colorbar style={
          width=.15cm, xshift=-0.9em, 
          point meta min=\cmin,point meta max=\cmax,
        },
        ]
        \addplot graphics [
        xmin=\xmin, xmax=\xmax,  ymin=\xmin, ymax=\xmax,
        ] {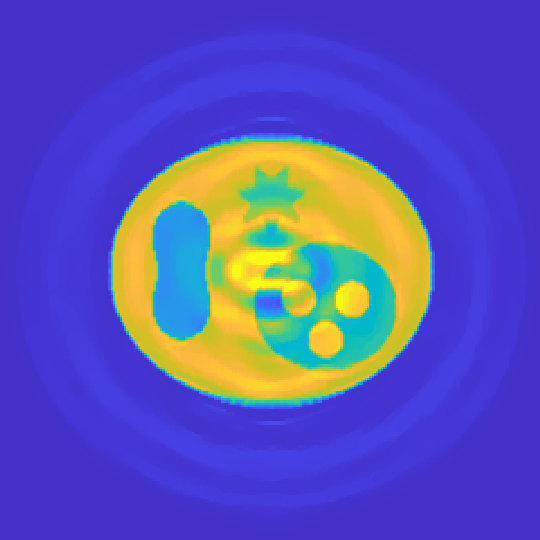};
      \end{axis}
    \end{tikzpicture}
    \caption{HIO/CG and TVd\\ $\lambda=0.02$, $J_{\mathrm{TV}}=20$\\PSNR 29.88, SSIM 0.713\\ 31\,sec}
  \end{subfigure}\hspace{-6pt}
  \begin{subfigure}[t]{.33\textwidth}\centering
    \begin{tikzpicture}
      \begin{axis}[
        width=\modelwidth,
        height=\modelwidth,
        enlargelimits=false,
        scale only axis,
        axis on top,
        colorbar,colorbar style={
          width=.15cm, xshift=-0.9em, 
          point meta min=\cmin,point meta max=\cmax,
        },
        ]
        \addplot graphics [
        xmin=\xmin, xmax=\xmax,  ymin=\xmin, ymax=\xmax,
        ] {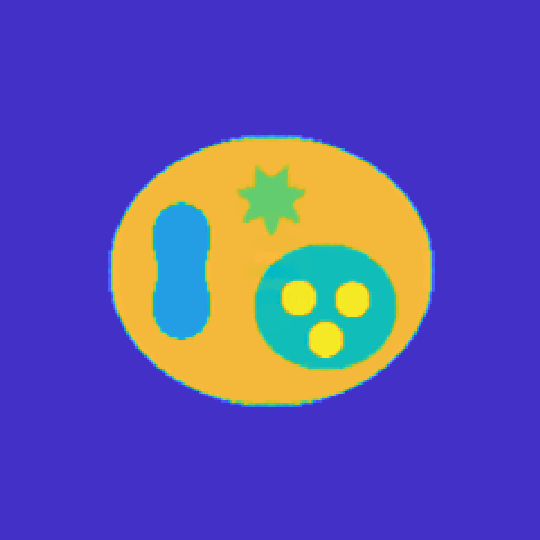};
      \end{axis}
    \end{tikzpicture}
    \caption{HIO/PD with TV\\
      \mbox{$\lambda{=}0.01$, $J_{\mathrm{IO}}{=}200$,
        $J_{\mathrm{PD}}{=}5$}
      \\PSNR 34.92, SSIM 0.994\\ 3\,h 48\,min}
  \end{subfigure}
  \caption{Slice $x_3=0.35$ of the 3D HIO phase retrieval using CG in Alg.~\ref{alg:cg} and PD
  with TV in Alg.~\ref{alg:pd-back}.  The HIO/CG result is
  TV-denoised.
  \label{fig:cell3d-pr}}
\end{figure}

\begin{figure}[t]\centering
\pgfmathsetmacro{\cmin} {-0.05} \pgfmathsetmacro{\cmax} {1.05}
\begin{subfigure}[t]{.33\textwidth}
  \begin{tikzpicture}
    \begin{axis}[
      width=\modelwidth,
      height=\modelwidth,
      enlargelimits=false,
      scale only axis,
      axis on top,
      colorbar,colorbar style={
        width=.15cm, xshift=-0.9em, 
        point meta min=\cmin,point meta max=\cmax,
      },
      ]
      \addplot graphics [
      xmin=\xmin, xmax=\xmax,  ymin=\xmin, ymax=\xmax,
      ] {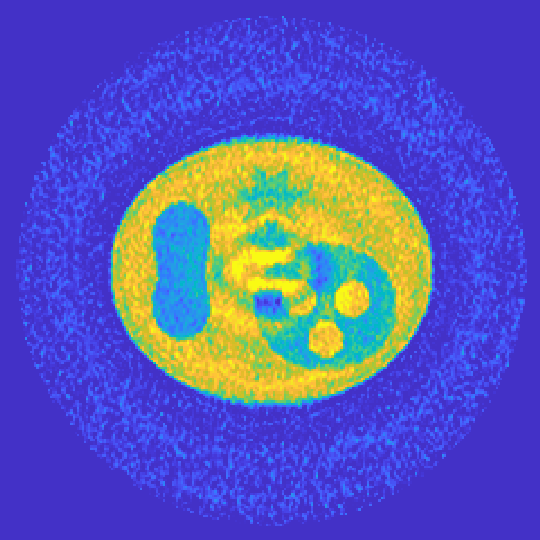};
    \end{axis}
  \end{tikzpicture}
  \caption{HIO/CG reconstruction\\  $J_{\mathrm{IO}}=10$, $J_{\mathrm{CG}}=5$\\PSNR 24.71, SSIM 0.589}
\end{subfigure}\hspace{-6pt}
\begin{subfigure}[t]{.33\textwidth}\centering
  \begin{tikzpicture}
    \begin{axis}[
      width=\modelwidth,
      height=\modelwidth,
      enlargelimits=false,
      scale only axis,
      axis on top,
      colorbar,colorbar style={
        width=.15cm, xshift=-0.9em, 
        point meta min=\cmin,point meta max=\cmax,
      },
      ]
      \addplot graphics [
      xmin=\xmin, xmax=\xmax,  ymin=\xmin, ymax=\xmax,
      ] {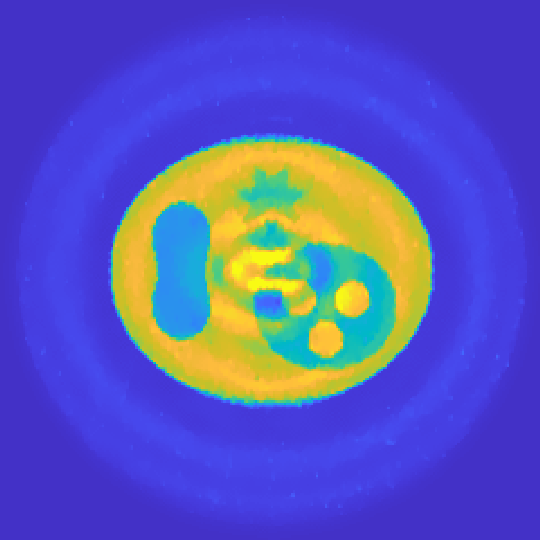};
    \end{axis}
  \end{tikzpicture}
  \caption{HIO/CG and TVd\\ $\lambda=0.05$, $J_{\mathrm{PD}}=20$\\PSNR 27.16, SSIM 0.623}
\end{subfigure}\hspace{-6pt}
\begin{subfigure}[t]{.33\textwidth}\centering
  \begin{tikzpicture}
    \begin{axis}[
      width=\modelwidth,
      height=\modelwidth,
      enlargelimits=false,
      scale only axis,
      axis on top,
      colorbar,colorbar style={
        width=.15cm, xshift=-0.9em, 
        point meta min=\cmin,point meta max=\cmax,
      },
      ]
      \addplot graphics [
      xmin=\xmin, xmax=\xmax,  ymin=\xmin, ymax=\xmax,
      ] {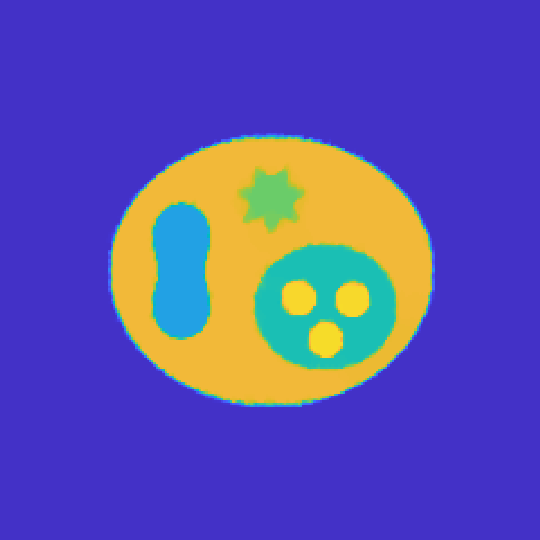};
    \end{axis}
  \end{tikzpicture}
  \caption{HIO/PD with TV\\
    \mbox{$\lambda{=}0.05$, $J_{\mathrm{IO}}{=}200$,
      $J_{\mathrm{PD}}{=}5$}
    \\PSNR 34.10, SSIM 0.993}
\end{subfigure}
\caption{Slice $x_3=0.35$ of the 3D HIO phase retrieval with 5\%
  Gaussian noise using CG in Alg.~\ref{alg:cg} and PD
  with TV in Alg.~\ref{alg:pd-back}.  The HIO/CG result is
  TV-denoised.
  \label{fig:cell3d-pr-noise}}
\end{figure}

For the three-dimensional simulations, we impose the support
constraint $\rs = 40$ in \eqref{eq:f-projection}.  The reconstructions
are shown in \autoref{fig:cell3d-pr}.  There are some artifacts in the
HIO/CG reconstruction.  Raising the number of HIO iterations in this
case promotes further noise artifacts, since the data have been
computed using a discretization of the convolution \eqref{eq:Dtot-op}.
The slight differences between the data and reconstruction model
avoids the inverse crime.  The HIO method with PD and TV overcomes
this issue and produces excellent reconstructions comparable to
the known-phase inversion.  Adding more noise to the data has a similar
effect to the reconstruction as in the two-dimensional case, see
\autoref{fig:cell3d-pr-noise}.

\subsection{Wave propagation--backpropagation method}
\label{sec:MD}

In contrast to applying a phase retrieval method to the full ODT model
$\zb D^\mathrm{tot}$, an alternative approach proposed by Maleki \&
Devaney \cite{MalDev93} takes into account the physical nature of
diffraction tomography.  We call it the MD approach.  It is also based
on the input-output methods of \autoref{sec:phase-retrieval-odt}, but instead of
solving the whole ODT problem in each step, only the wave propagation
in free space is considered.  In particular, we define the free space
wave propagation operator $\mathcal P$ for $v\colon \R^{d-1}\to\C$ by
its Fourier transform
\begin{equation*}
\mathcal F [\mathcal Pv](\by')
= \e^{\i\, \kappa(\by')\, \rM}\, \mathcal F[ v](\by')
,\qquad \by'\in \R^{d-1},
\end{equation*}
where $\mathcal F$ denotes the $d-1$ dimensional Fourier transform.
The inverse, the so-called free space backpropagation
$\mathcal P^{-1}$, applied to the field $\us_t(\cdot,\rM)$ is
assumed to approximately match $\us_t(\cdot,0)$, which holds true
for very thin samples.  Then the support of
$$
\us_t(\cdot,0)
\approx \mathcal P^{-1} \big( \utot_t(\cdot,\rM) - \ui(\cdot,\rM) \big)
$$
is assumed to be contained in the ball of radius $\rs$.

Overall, the MD approach consists of two steps.  First, we apply
\autoref{alg:err-red} with $\bm D$ replaced by
$\mathcal P + \e^{\i k_0 \rM}$ to the data $\abs{\utot_t(\cdot,\rM)}$
and impose that $\us(\cdot,0)$ has a compact support of radius $\rs$.
In contrast to previous sections, we can assume neither the
non-negativity nor the realness of $\us(\cdot,0)$.  In the second
step, we use an inversion method of \autoref{sec:indft} to recover
$f$ using the retrieved phase of $\utot$ from the first step.

As pointed out in \cite{MalDev93}, this method works only well if the
distance $\rM$ of the measurement plane is very large, for instance
$\rM>1000$.  To obtain decent reconstructions, we require larger
measurement plane sizes $\lM$.  Consequently $N$ has to becomes
larger if the distance between neighboring grid points should stay the same as
before.

\begin{figure}[tb]\centering
  \pgfmathsetmacro{\cmin} {-0.05} \pgfmathsetmacro{\cmax} {0.55}
  \begin{subfigure}[t]{.25\textwidth}
    \begin{tikzpicture}
      \begin{axis}[
        width=.8*\modelwidth,
        height=.8*\modelwidth,
        enlargelimits=false,
        scale only axis,
        axis on top,
        ]
        \addplot graphics [
        xmin=\xmin, xmax=\xmax,  ymin=\xmin, ymax=\xmax,
        ] {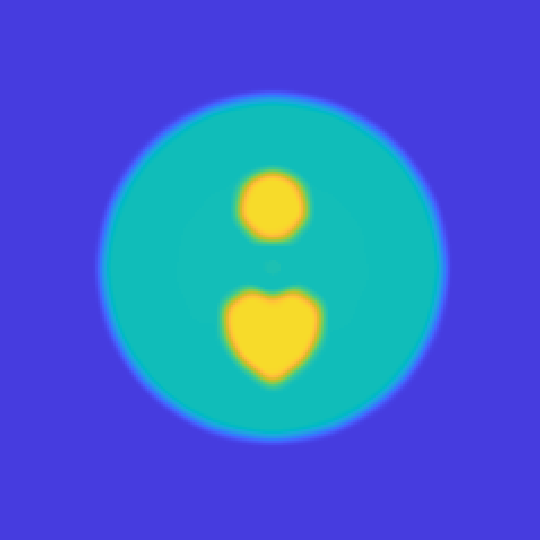};
      \end{axis}
    \end{tikzpicture}
    \caption{HIO/PD with TV\\$J_{\mathrm{IO}}=50$, $J_{\mathrm{CG}}=5$\\ PSNR~33.81\\ SSIM~0.945\\
    34\,sec}
  \end{subfigure}\hfill
  \begin{subfigure}[t]{.25\textwidth}\centering
    \begin{tikzpicture}
      \begin{axis}[
        width=.8*\modelwidth,
        height=.8*\modelwidth,
        enlargelimits=false,
        scale only axis,
        axis on top,
        ]
        \addplot graphics [
        xmin=\xmin, xmax=\xmax,  ymin=\xmin, ymax=\xmax,
        ] {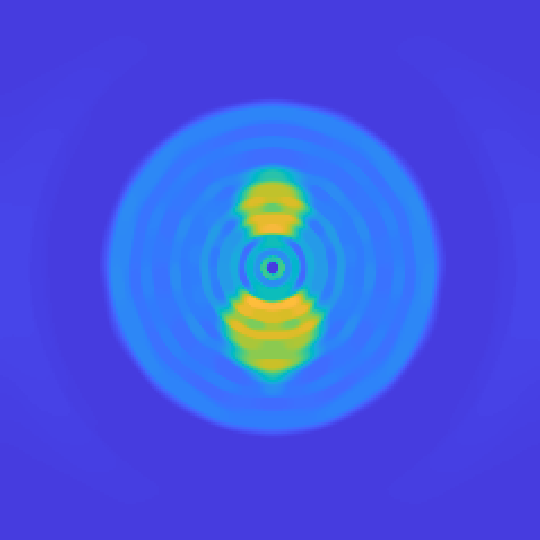};
      \end{axis}
    \end{tikzpicture}
    \caption{MD retrieval\\ $J_{\mathrm{IO}} = 2\cdot10^4$\\ PSNR~21.74\\ SSIM~0.771\\
    3\,min~16\,sec}
  \end{subfigure}\hfill
  \begin{subfigure}[t]{.25\textwidth}\centering
    \begin{tikzpicture}
      \begin{axis}[
        width=.8*\modelwidth,
        height=.8*\modelwidth,
        enlargelimits=false,
        scale only axis,
        axis on top,
        ]
        \addplot graphics [
        xmin=\xmin, xmax=\xmax,  ymin=\xmin, ymax=\xmax,
        ] {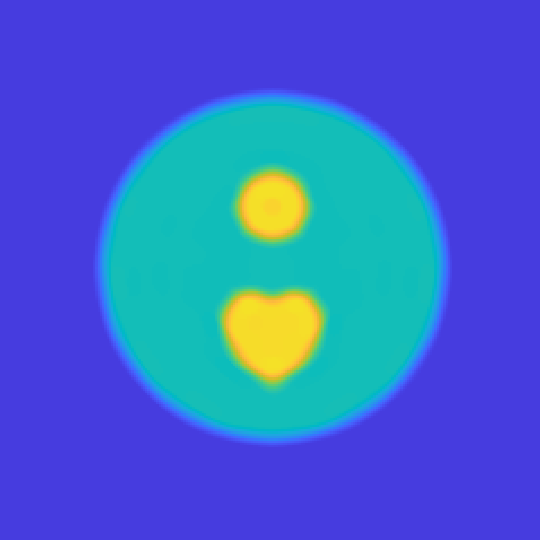};
      \end{axis}
    \end{tikzpicture}
    \caption{MD retrieval\\ $J_{\mathrm{IO}} = 5\cdot10^4$\\ PSNR~33.37\\ SSIM~0.937\\ 8\,min}
  \end{subfigure}\hfill
  \begin{subfigure}[t]{.25\textwidth}\centering
  \begin{tikzpicture}
    \begin{axis}[
      width=.8*\modelwidth,
      height=.8*\modelwidth,
      enlargelimits=false,
      scale only axis,
      axis on top,
      ]
      \addplot graphics [
      xmin=\xmin, xmax=\xmax,  ymin=\xmin, ymax=\xmax,
      ] {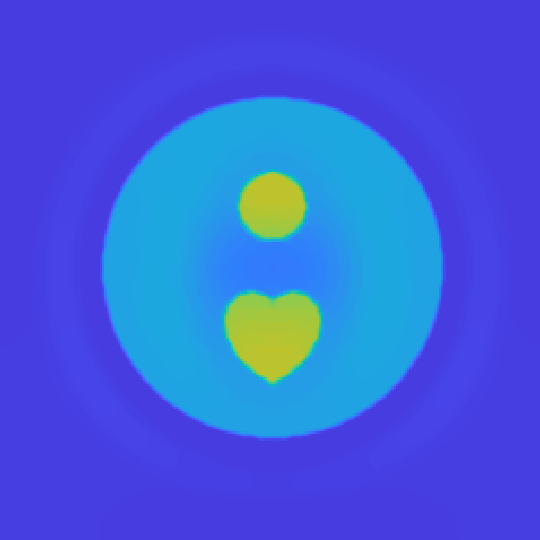};
    \end{axis}
  \end{tikzpicture}
  \caption{MD retrieval\\ $J_{\mathrm{IO}} = 5\cdot10^6$\\ PSNR~25.78\\ SSIM~0.867\\ 13\,h \label{fig:heart1-pr-md-rM400}}
  \end{subfigure}
  \caption{MD phase retrieval with $\rM=1000$ (a)-(c) and $\rM=400$ (d).
    \label{fig:heart1-pr-md}}
\end{figure}

For our numerical simulations, we chose $\lM=240$ and $\rM=1000$.
Further, we raise the number of grid points to $N=960$.  The rest of
the setting remains the same as in \autoref{sec:full-diffr-model}.
The second step of the MD method is conducted with the inversion via
PD with TV using $\lambda= 0.05$ and $J_{\mathrm{PD}} = 100$.  The
results are shown in \autoref{fig:heart1-pr-md}.  We see that the MD
approach requires a very high number of iterations to find a good
solution.  Even 20\,000 iterations are not sufficient, but 50\,000
yield a reasonable reconstruction.

For a smaller distance $\rM$, even more iterations are necessary.
With $\rM=400$, we stopped after 5 million iterations, which took
about 13 hours, and the reconstruction showed the right shape, but the
level was still too low, see \autoref{fig:heart1-pr-md-rM400}.  In
case of noisy data, we could not achieve a reasonable reconstruction
with the MD approach.

\section{Conclusion}
\label{sec:conclusion}

We have compared reconstruction algorithms in ODT based on the Born
approximation.  For known-phase measurements, we compared the
reconstruction based on three algorithms: the discrete
backpropagation, the CG method, and the PD method with TV regularization.  The
backpropagation is the fastest, but the image quality was inferior
compared to the others.  The CG method was better and could still be
improved when combined with TV denoising.  In this case, the results
are similar to the PD algorithm, which is usually the most accurate but also the
slowest method.

For unknown phase, which causes the reconstruction problem to be more
ill-posed, we achieved good results by applying the HIO algorithm,
where the inner reconstruction is performed with the PD method.  The
image quality is comparable to the full phase inversion.  The
implementation with the CG method falls behind in terms of quality
but, as it is much faster, gives good initial solutions for the HIO/PD
method.  By this procedure, the time requirements of PD are
significantly reduced.

\subsection*{Acknowledgments}

R.\ Beinert greatfully acknowledges funding by the BMBF under the project ``VI-Screen'' (13N15754).
M.\ Quellmalz greatfully acknowledges funding by the DFG (STE 571/19-1, project number 495365311)
within SFB F68 (“Tomography across the Scales”).

\printbibliography[heading=bibintoc]

\end{document}